\documentclass[11pt]{article}

\usepackage{amssymb}
\usepackage{amsmath}
\usepackage{amsthm}
\usepackage{graphicx}

\newcommand{\R}{\mathbb R}
\newcommand{\E}{\mathbb E}
\newcommand{\N}{\mathbb N}
\newcommand{\PO}{\mathbb P}

\newcommand{\cI}{{\cal I}}
\newcommand{\cP}{{\cal P}}
\newcommand{\bfe}{{\bf e}}
\newcommand{\bfi}{{\bf i}}

\newcommand{\ds}{\displaystyle}
\newcommand{\qedwhite}{\hfill \ensuremath{\Box}}

\newtheorem{theorem}{Theorem}[section]
\newtheorem{proposition}{Proposition}[section]

\usepackage{xcolor,todonotes}
\definecolor{otherblue}{rgb}{0,0.3,0.6}

\def\TolX{{\eta_{X}}}
\def\TolXm1{{\eta_{X_{-1}}}}
\def\TolXzero{{\eta_{X_0}}}
\def\TolXone{{\eta_{X_1}}}
\def\TolXtwo{{\eta_{X_2}}}
\def\TolXkm1{{\eta_{X_{k-1}}}}
\def\TolXk{{\eta_{X_k}}}
\def\TolXK{{\eta_{X_K}}}
\def\TolY{{\eta_{Y}}}
\def\TolYm1{{\eta_{Y_{-1}}}}
\def\TolYzero{{\eta_{Y_0}}}
\def\TolYone{{\eta_{Y_1}}}
\def\TolYtwo{{\eta_{Y_2}}}
\def\TolYkm1{{\eta_{Y_{k-1}}}}
\def\TolYk{{\eta_{Y_k}}}
\def\TolYK{{\eta_{Y_K}}}
\def\TolYKmk{{\eta_{Y_{K-k}}}}

\def\matlab{{\sc Matlab}}
\def\vfld#1{{#1}}
\def\bfy{y}
\def\ML{\mathrm{ML}}
\def\SL{\mathrm{SL}}

\author{J.~Lang, R.~Scheichl and D.~Silvester}
\title{A fully adaptive multilevel stochastic collocation strategy for
solving elliptic PDEs with random data}
\author{
Jens Lang\footnote{Corresponding author.}\\
{\small \it Technische Universit\"at Darmstadt, Department of Mathematics} \\
{\small \it Dolivostra{\ss}e 15, 64293 Darmstadt, Germany}\\[2mm]
Robert Scheichl\\
{\small \it Ruprecht-Karls-Universit\"at Heidelberg}\\
{\small \it Institute for Applied Mathematics and Interdisciplinary Center for Scientific Computing}\\
{\small \it Im Neuenheimer Feld 205, 69120 Heidelberg, Germany}\\[2mm]
David Silvester\\
{\small \it The University of Manchester, Department of Mathematics}\\
{\small \it Manchester M13 9PL, United Kingdom}
}
\date{December 4, 2019}
\begin{document}
\maketitle

\begin{abstract}
We propose and analyse a fully adaptive strategy
for solving elliptic PDEs with random data in this work.
A hierarchical sequence of adaptive mesh refinements for the spatial
approximation is combined with adaptive anisotropic sparse Smolyak
grids in the stochastic space in such a way as to minimize the computational
cost. The novel aspect of our strategy  is that the hierarchy of spatial approximations
is sample dependent so that the computational effort at each collocation point can be
optimised individually. We outline a rigorous analysis for the convergence and
computational complexity of the adaptive multilevel algorithm and we provide optimal choices for error tolerances at each level. Two numerical examples demonstrate the reliability of the error control and the significant decrease in the complexity that arises when compared
to single level algorithms and multilevel algorithms that employ adaptivity solely in the spatial discretisation or in the collocation procedure.
\end{abstract}



\section{Introduction}
A central task in computational science and engineering is the
efficient numerical treatment of partial differential equations (PDEs)
with uncertain input data (coefficients, source term, geometry,
etc). This has been a very active area of research
in recent years with a host of competing ideas.  Our focus here is on
elliptic PDEs with
random input data in a standard (pointwise) Hilbert space setting: Find
 $u(\cdot,y)\in V$ such that
\begin{align}
\label{eqn_rand1}
 a_y( u(\cdot,y), v)   &= \ell_y(v) \quad \forall v\in V, \,y \in \Gamma
\end{align}
where $ a_y(\cdot, \cdot): V \times V \to \R$ is a parameter-dependent
inner product in a Hilbert space $V$ and  $\ell_y(v): V \to \R$ is a
parameter-dependent bounded linear functional.
In general, the solution will be represented as $u(x,y):D\times\Gamma\rightarrow\R$
where $D\subset\R^d$, $d=1,2,3$, is the deterministic, bounded (physical) domain
and $\Gamma = \Gamma_1 \times \Gamma_1 \times \cdots \Gamma_N$
is a (stochastic) parameter space of finite dimension $N$ (finite noise assumption).
The component parameters  $y_1,\ldots,y_n$  will be associated with independent
random variables  that have  a joint probability density function
$\rho(y)=\Pi_{n=1}^N \hat{\rho}(y_n) \in L^\infty(\Gamma)$ such that
$\hat{\rho}_n:[-1,1]\rightarrow\R$.

The most commonly used approach to solve \eqref{eqn_rand1} is the
Monte Carlo (MC) method based on random samples $y^{(m)} \in
\Gamma$. It has a dimension-independent convergence rate and leads to
a natural decoupling of the stochastic and spatial dependency. While  MC methods  are well
suited to adaptive spatial refinement,  they are not able to exploit any smoothness or special structure in the parameter dependence.
The standard, single-level stochastic collocation method (see, for example, \cite{BabuskaNobileTempone2007,NobileTemponeWebster2008b,
NobileTemponeWebster2008a,XiuHesthaven2005})
for \eqref{eqn_rand1} is similar to the MC method in that it involves independent,
finite-dimensional spatial approximations $u_h(y^{(m)}) \approx
u(y^{(m)})$ at a set $\{y^{(m)}\}_{m=1,\ldots,M}$ of deterministic
sampling points in $\Gamma$, {so as to} construct an interpolant
\begin{equation}
u_{M,h}^{(\SL)}(x,y) = \cI[u_h](x,y) = \sum_{m=1}^{M} u_m(x)\phi_m(y),
\end{equation}
in the polynomial space $\cP_M=\mbox{span}\{\phi_m\}_{m=1,\ldots,M}
\subset L^2_\rho(\Gamma)$ with
basis functions $\{\phi_m\}_{m=1,\ldots,M}$. The coefficients $u_m(x)$ are
determined by the interpolating condition $\cI[u_h](x,y^{(m)})=u_h(x,y^{(m)})$, for
$m=1,\ldots,M$. The quality of the interpolation process depends on the
accuracy of the spatial approximations $u_h(y^{(m)})$ and on the number
of collocation points $M$, which typically, grows rapidly with increasing stochastic
dimension $N$.

Two concepts in the efficient treatment of \eqref{eqn_rand1}
that have shown a lot of promise and are (to some extent)
complementary are (i) adaptive sparse-grid interpolation methods that
are able to exploit smoothness with respect to the uncertain
parameters and (ii) multilevel
approaches that aim at reducing the computational cost through
a hierarchy of spatial approximations.
Adaptive sparse-grid methods, where the set of sample points
is adaptively generated, can be traced back to Gerstner \&
Griebel~\cite{GerstnerGriebel2003} and have been extensively tested in
collocation form; see~\cite{AlmeidaOden2010,
ChkifaCohenSchwab2014,GuignardNobile2018,NobileTamelliniTeseiTempone2016,SchiecheLang2014}.
We note that parametric adaptivity has also been explored in  a Galerkin framework; see \cite{BespalovPowellSilvester2014,BespalovSilvester2016,CohenChifkaSchwabDevore2013, EigelGittelssonSchwabZander2014,EigelGittelssonSchwabZander2015},
and that there are a number of recent papers aimed at proving
dimension-independent convergence; see \cite{BachmayrCohenDungSchwab2017,
BespalovPraetoriusRocchiRuggeri2019,NobileTamelliniTempone2016,Zech2018,ZechDungSchwab2018}.

The multilevel Monte Carlo (MLMC) method was originally proposed as an
abstract variance reduction technique by Heinrich~\cite{Heinrich2001}, and independently
for stochastic differential equations in mathematical finance by Giles~\cite{Giles2008}.
In the context of uncertainty quantification it was developed  by
Barth et al.~\cite{BarthSchwabZollinger2011} and Cliffe et al.~\cite{CliffeGilesScheichlTeckentrup2011} and
extended  to stochastic collocation sampling methods soon thereafter
 by Teckentrup et al.~\cite{TeckentrupJantschWebsterGunzburger2012}.
The methodology in this paper can be viewed as an extension of this body of work.

Problems where the PDE solution develops singularities (or at least
steep gradients) in random locations in the spatial domain are a challenging
feature in many real-world applications. Spatial adaptivity driven by  a posteriori error estimates has been investigated within single-level stochastic collocation in
\cite{SchiecheLang2014} and has  been considered in a MLMC framework in~\cite{EigelMerdonNeumann2016}.
To date, however, MC-based methods  \cite{DetommasoDodwellScheichl2019,ElfversonHellmannMalqvist2016,KornhuberYouett2018,Youett2018}
are the only methods  where  the spatial approximation is locally refined at each
parametric sample point.
In what follows, we will combine a pointwise adaptive mesh refinement for the approximation
of the spatial approximations $u_h(y^{(m)})$ with an adaptive (anisotropic) sparse Smolyak grid,
in order to improve the efficiency of the multilevel method to reach a user-prescribed
tolerance for the accuracy of the multilevel interpolant, as well as
for the computation  of associated quantities of interest.
By rigorous control of both error components,
the fully adaptive algorithm proposed herein
 is able to exploit smoothness and structure in the parameter
dependence while  also exploiting fully the computational gains due to the spatially adapted
hierarchy of spatial approximations at each collocation point. Under
certain assumptions on the efficiency of the two-component adaptive schemes, the
complexity of the new multilevel algorithm can be estimated rigorously. Two
specifically chosen numerical examples will be looked at in the final section
to demonstrate its efficiency.

\section{Methodology}
\label{methodology}

Following the paper of
Teckentrup et al.~\cite{TeckentrupJantschWebsterGunzburger2012},
we consider steady-state diffusion problems with uncertain parameters.
Thus, we will assume that $V := H^1_0(D)$ and that the variational
formulation (\ref{eqn_rand1})  admits a unique solution in the weighted
Bochner space
\begin{equation}
L_\rho^2(\Gamma;H^1_0(D))=\{v:\Gamma\rightarrow H^1_0(D)
\mbox{ measurable: }
\int_\Gamma \|v(y,\cdot) \|^2_{H^1_0(D)} \, \rho(y) \, dy < \infty \}
\end{equation}
with corresponding norm
\begin{equation}
\| v\|^2_{L_\rho^2(\Gamma;H^1_0(D))} =
\int_\Gamma \|v(y,\cdot) \|^2_{H^1_0(D)}\, \rho(y) \, dy
=:\E\left[ \|v(y,\cdot) \|^2_{H^1_0(D)}\right].
\end{equation}
The structure  of our adaptive algorithm is
identical to that introduced in~\cite{TeckentrupJantschWebsterGunzburger2012}.
The distinctive aspects are discussed in this section.

\subsection{Adaptive spatial approximation}
We will focus on finite element approximation in the spatial domain $D$ using classical adaptive
methodology.  That is, we  execute the loop
\begin{equation} \label{eq:modules}
    \mathsf{SOLVE}\ \Longrightarrow\ \mathsf{ESTIMATE}\ \Longrightarrow\ \mathsf{MARK}\ \Longrightarrow\ \mathsf{REFINE}
 \end{equation}
until the estimate of the error in the second step  is less than a prescribed tolerance.
Next, we let $\{\TolXk\}_{k=0,...,K}$ be a decreasing sequence of tolerances with
\begin{equation}
\label{seq-tolx}
1 \ge \TolXzero > \TolXone > \cdots > \TolXk > \cdots > \TolXK > 0.
\end{equation}
Then, for each fixed parameter $y\in\Gamma$, we run our adaptive
strategy to compute approximate spatial solutions $u_k(y)\in V_k(y)$ on
a sequence of nested subspaces
\begin{equation}
V_0(y)\subset V_1(y)\subset\cdots\subset V_K(y)\subset H^1_0(D).
\end{equation}
This is the point of departure from  the previous work~\cite{TeckentrupJantschWebsterGunzburger2012},
wherein  the sequence of approximation spaces $V_0\subset V_1 \subset\cdots\subset V_K$
could be chosen in a general way, but in a manner that was fixed for all
collocation points $y\in \Gamma$.

Starting from the pointwise error estimate
\begin{equation}
\label{equ-space-err-ptwse}
\| u(\cdot,y)- u_k(y) \|_{H^1_0(D)} \le C_k(y)\cdot \TolXk,\quad
k=0,\ldots,K,
\end{equation}
and supposing measurability of the discrete spaces $V_k(y)$ and bounded
second moments of $C_k(y)$, we directly get the following error bound
\begin{equation}
\label{equ-space-err}
\| u - u_k \|_{L^2_\rho(\Gamma,H^1_0(D))} \le C_X\cdot \TolXk,\quad
k=0,\ldots,K,
\end{equation}
with a constant
\begin{equation}
\label{def_cx}
C_X:=\max_{k=0,\ldots,K} \left( \int_\Gamma C^2_k(y) \, \rho(y) \, dy \right)^{1/2}
\end{equation}
that does not depend on $y$ and $k$.
Adaptive algorithms proposed by D\"orfler~\cite{Doerfler1996} and
Kreuzer\cite{KreuzerSiebert2011} converge for fixed $y\in\Gamma$ and $\TolXk\rightarrow 0$.
The constant $C_X$ in (\ref{def_cx}) is related to the effectivity of the a
posteriori error estimation strategy.
Values close to one can be obtained using hierarchical error estimators (see, for example,
Bespalov et al.~\cite{BespalovRocchiSilvester2019}) or using gradient recovery techniques in the asymptotic regime.

\subsection{Adaptive stochastic interpolation}
Let us assume $u\in C^0(\Gamma;H^1_0(D))$ and denote by
$\{\cI_{M_k}\}_{k=0,1,\ldots}$ a sequence of interpolation operators
\begin{equation}
\label{stoch-intpl-op}
\cI_{M_k}: C^0(\Gamma) \rightarrow L^2_\rho(\Gamma)
\end{equation}
with $M_k$ points from the $N$-dimensional space
$\Gamma$. We construct each of these
operators by a hierarchical sequence of
one-dimensional Lagrange interpolation operators with the anisotropic
Smolyak algorithm, which was introduced by Gerstner \& Griebel~\cite{GerstnerGriebel2003}.
The method is dimension adaptive, using
the individual surplus spaces in the multi-dimensional hierarchy as
natural error indicators.

Let $\{\TolYk \}_{k=0,...,K}$ be a second sequence of tolerances,
a priori not necessarily decreasing.
Under suitable regularity assumptions for the uncertain data
(see e.g. Babu{\v s}ka et al.~\cite[Lem.~3.1, 3.2]{BabuskaNobileTempone2007}),
we can assume that there exist numbers $M_k$, $k=0,1,\ldots,K$, and a
constant $C_Y > 0$ not depending on $k$ such that
\begin{equation}
\label{equ-stoerr}
\| (u_k-u_{k-1}) - \cI_{M_{K-k}} [u_k-u_{k-1}]\|_{L^2_\rho(\Gamma;H^1_0(D))}
\le C_Y\cdot \TolYKmk,\quad k=0,\ldots,K,
\end{equation}
where, for simplicity, we set $u_{-1}=0$.

Since \eqref{equ-space-err-ptwse} implies
$\|u_k-u_{k-1}\|_{L^2_\rho(\Gamma;H^1_0(D))}\le C\cdot
\TolXkm1\,$, which is
decreasing as $k\rightarrow\infty$, we can expect that for higher
$k$ it suffices to use less accurate
interpolation operators, i.e., smaller numbers $M_{K-k}$,
to achieve a required accuracy. Indeed this is the main motivation to set
up a multilevel interpolation approximation. Note that in this way, the tolerances
$\TolYKmk$ are strongly linked to the spatial tolerances $\TolXk$. We will make
this connection more precise and give suitable values for the
sequence of tolerances $\{\TolYk\}_{k=0,\ldots,K}$ in the next section.

\subsection{An adaptive multilevel strategy}
Given the sequences $\{u_k\}$ and $\{\cI_{M_k}\}$, we define
the multilevel interpolation approximation in the usual way
by
\begin{equation}
\label{mlipol}
u_K^{(\ML)} = \sum_{k=0}^{K} \cI_{M_{K-k}} [u_k - u_{k-1}] =
\sum_{k=0}^{K} \left( u_{M_{K-k},k}^{(\SL)} - u_{M_{K-k},k-1}^{(\SL)}\right).
\end{equation}
Observe that the most accurate interpolation operator $\cI_{M_K}$ is used
on the coarsest spatial approximation $u_0$ whereas the least accurate
interpolation operator $\cI_{M_0}$ is applied to the difference of the
finest spatial approximations $u_{K}-u_{K-1}$. The
close relationship between the spatial and stochastic
approximations  at index $k$ is  clearly visible in (\ref{mlipol}).
We note that the use of adaptive interpolation operators is also mentioned
in~\cite{TeckentrupJantschWebsterGunzburger2012}. The key difference
is that the hierarchy in~\cite{TeckentrupJantschWebsterGunzburger2012} is
based on the number of collocation points $M_k$ rather than the
tolerances $\{\TolYk\}$ employed here.

To show the convergence of the multilevel approximation $u_K^{(\ML)}$ to
the true solution $u$, we split the error into the sum of a spatial
discretization error and a stochastic interpolation error using the triangle inequality
\begin{equation}
\label{est-totalerr}
\| u - u_K^{(\ML)} \|_{L^2_\rho(\Gamma;H^1_0(D))} \le
\| u - u_K \|_{L^2_\rho(\Gamma;H^1_0(D))} +
\| u_K - u_K^{(\ML)}\|_{L^2_\rho(\Gamma;H^1_0(D))}.
\end{equation}
Due to (\ref{equ-space-err}) the first term on the right hand side
of \eqref{est-totalerr} is bounded
by $C_X\cdot \TolXK$. The aim is now to choose the tolerances $\{\TolYk\}$
in an appropriate way to reach the same accuracy for the second term.
From (\ref{equ-stoerr}), we
estimate the stochastic interpolation error as follows:
\begin{equation}
\label{est-stocherr}
\begin{array}{rll}
\| u_K - u_K^{(\ML)}\|_{L^2_\rho(\Gamma;H^1_0(D))} &\!=\!&
\ds\left\| \sum_{k=0}^{K} (u_k - u_{k-1}) -
\cI_{M_{K-k}}[u_k - u_{k-1}] \right\|_{L^2_\rho(\Gamma;H^1_0(D))} \\[5mm]
&& \hspace{-4cm}\le\ds\sum_{k=0}^{K} \left\| (u_k - u_{k-1}) -
\cI_{M_{K-k}}[u_k - u_{k-1}]\right\|_{L^2_\rho(\Gamma;H^1_0(D))}
\le \sum_{k=0}^{K} C_Y\cdot \TolYKmk.
\end{array}
\end{equation}
To obtain an accuracy of the same size as the spatial discretization
error, a first choice for the tolerances would be to simply demand
$\TolYk\le C_X\cdot \TolXK/((K+1)C_Y)$, for all
$k=0,\ldots,K$, and conclude that
\begin{equation}
\| u - u_K^{(\ML)} \|_{L^2_\rho(\Gamma;H^1_0(D))} \le 2\,C_X\cdot \TolXK,
\end{equation}
so that  the adaptive multilevel method converges for
$\TolXK\rightarrow 0$. However, the values for $\TolYk$ can be optimized
by minimizing the computational costs while keeping the desired accuracy.

At this point it is appropriate to consider the
computational cost, $C_\epsilon^{(\ML)}$, of the multilevel
stochastic collocation estimator $u_K^{(\ML)}$  required to achieve
an accuracy $\epsilon$. Thus, in order to quantify the contributions from the spatial discretization and the stochastic collocation, we will need to make two assumptions to
link the cost with the error bounds in (\ref{equ-space-err}) and
(\ref{equ-stoerr}). Let $A_k$ denote an upper bound for the cost to
solve the deterministic PDE
at sample point $y \in \Gamma$ with accuracy $\TolXk$. Then,
we assume that for all $k=0,\ldots,K$,
\begin{equation*}
\begin{array}{lrll}
\text{(A1)} \qquad & \quad A_k &\le& C_A \cdot \TolXk\!\!{}^{-s},\qquad\\[2mm]
\text{(A2)} \qquad & \quad C_Y \cdot \TolYKmk &=&
C_I(N)\,M_{K-k}^{-\mu}\,\TolXkm1 \qquad
\end{array}
\end{equation*}
as well as the special case
$$
\|u_0\|_{L^2_\rho(\Gamma;H^1_0(D))} \ \le \ \TolXm1:=const.
$$
Here, the constants $C_A>0$, $C_I(N)>0$ are independent of $y$, $k$, and
the rates $s,\mu>0$.
Note that we could also consider the exact cost per sample and introduce a sample dependent
constant $C_A(y)$ in (A1), but that would complicate the subsequent analysis.

Assumption (A1) usually holds for first-order adaptive spatial discretization methods
with $s=d$, when coupled with optimal linear solvers such as multigrid. The
factors on the right-hand side in (A2) reflect best the
convergence of the sparse grid approximations in (\ref{equ-stoerr}) with respect to the total number $M_{K-k}$ of collocation points,
see Nobile et al.~\cite{NobileTemponeWebster2008b,NobileTemponeWebster2008a}
or \cite[Theorem 5.5]{TeckentrupJantschWebsterGunzburger2012}. To estimate the
difference $u_k-u_{k-1}$, we use the
fact that $\|u_k-u_{k-1}\|_{L^2_\rho(\Gamma;H^1_0(D))}\le C\cdot \TolXkm1$
with a constant $C>0$ close to $C_X$. It follows from
$\|u-u_{k-1}\|_{H^1_0(D)}\approx \|u_k-u_{k-1}\|_{H^1_0(D)}$, which is the basis for the very good performance of hierarchical error estimators. We will absorb $C$ into $C_I$ in the sequel.

The rate $\mu$ depends in general on the dimension $N$. Theoretical
results for the anisotropic classical Smolyak algorithm are given in
\cite[Thm.~3.8]{NobileTemponeWebster2008b}. However, it has recently
been shown in \cite{Zech2018,ZechDungSchwab2018} that under certain
smoothness assumptions $\mu$ can be independent of $N$.

An upper bound for the total computational cost
of the approximation $u_K^{(\ML)}$ can then be defined as
\begin{equation}
\label{mlcost}
\ds C^{(\ML)} \,\leq\, \sum_{k=0}^{K} M_{K-k}\,(A_k+A_{k-1}),
\end{equation}
with $A_{-1}:=0$.
In a first step, we will consider a general sequence $\{\TolXk\}_{k=0,\ldots,K}$ without
defining a decay rate a priori.
Following the argument in~\cite{TeckentrupJantschWebsterGunzburger2012},
with a priori estimates for the spatial error  replaced by  tolerances, then leads to an estimate of $\epsilon$-cost $C_\epsilon^{(\ML)}$ and identifies
a set of optimal tolerances $\TolYk$ in (\ref{equ-stoerr}).
\begin{theorem}
Let $\TolXzero, \TolXone, \ldots,$ be a decreasing sequence of spatial tolerances satisfying (\ref{seq-tolx}). Suppose assumptions (A1)
and (A2) hold. Then, for any $\epsilon$, there exist an integer
$K = K(\epsilon)$ and a sequence of tolerances
$\{\TolYk\}_{k=0,\ldots,K}$ in (\ref{equ-stoerr}) such that
\begin{equation}
\| u - u_K^{(\ML)}\|_{L^2_\rho(\Gamma;H^1_0(D))} \le \epsilon
\end{equation}
and
\begin{equation}
\label{epscost}
\ds C_\epsilon^{(\ML)} \le C\cdot \big(G_K(\mu)\big)^{\frac{\mu+1}{\mu}}\,
\epsilon^{-\frac{1}{\mu}}
+ C_A\, \sum_{k=0}^{K}F_k(s) \, \TolXkm1
\end{equation}
with $C=C_A\,(2\,C_I)^{\frac{1}{\mu}}$ and
\begin{align}
\label{Fkdef}
F_k(s) &= \left( \TolXk\!\!{}^{-s}+ \TolXkm1\!\!\!\!\!\!\!{}^{-s}\right)\,\TolXkm1\!\!\!\!\!\!{}^{-1},
\quad k=0,\ldots,K,\\
\label{Gkdef}
G_K(\mu ) &= \sum_{k=0}^{K} \big(F_k(s)\big)^{\frac{\mu}{\mu+1}}\,\TolXkm1.
\end{align}
where for ease of notation we set $\TolXm1=\infty$.
The optimal choice for the tolerances $\TolYk$ is then given by
\begin{equation}
\label{opt-toly}
\TolYKmk = \big(2\,C_Y\,G_K(\mu)\big)^{-1} \big(F_k(s)\big)^{\frac{\mu}{\mu+1}}\,
\TolXkm1\,\epsilon.
\end{equation}
The utility of (\ref{opt-toly}) is that   near-optimal tolerances can be readily computed if estimates
 of the constants and the rates in (A1) and (A2) are available.
\end{theorem}
\noindent {\it Proof:} As in the convergence analysis above, we split the error
and make sure that both the spatial discretization error and the stochastic
interpolation error are bounded by $\epsilon/2$. First, we choose an
appropriate $K\ge 0$ and $\TolXk$ such that $C_X\cdot \TolXk<\epsilon/2$. This is,
of course, always possible and fixes the number $K=K(\epsilon)$ as a
function of $\epsilon$. Next we determine
the set $\{M_k\}_{k=0,\ldots,K}$ so that the computational cost in
(\ref{mlcost}) is minimized subject to the requirement that the
stochastic interpolation error is bounded by $\epsilon/2$. Using assumptions
(A1) and (A2), this reads
\begin{equation}
\ds\min_{M_0,\ldots,M_K} \;\sum_{k=0}^{K} C_A\cdot M_{K-k}
\left( \TolXk\!\!{}^{-s}+ \TolXkm1\!\!\!\!\!\!\!{}^{-s}\right)
\quad s.t.\quad \sum_{k=0}^{K} C_I\cdot M_{K-k}^{-\mu}\,\TolXkm1
=\frac{\epsilon}{2}.
\end{equation}
Application of the Lagrange multiplier method with all $M_k$ treated as
continuous variables as in Giles~\cite{Giles2008} gives the optimal choice for the number of samples
\begin{equation}
M_{K-k} = \big(2\,C_I\,G_K(\mu)\big)^{\frac{1}{\mu}}
\big(F_k(s)\big)^{-\frac{1}{\mu+1}}\,\epsilon^{-\frac{1}{\mu}}
\end{equation}
with $F_k$ and $G_K$ defined in (\ref{Fkdef}) and  (\ref{Gkdef}).
To ensure that $M_{K-k}$ is an integer, we round up to the next integer.
The $\epsilon$--cost of the multilevel approximation can then be estimated as follows
\begin{align*}
C_\epsilon^{(\ML)} &\le \,  \sum_{k=0}^{K} C_A\cdot (M_{K-k}+1)
\left( \TolXk\!\!{}^{-s}+  \TolXkm1\!\!\!\!\!\!\!{}^{-s}\right)\\
&=\, \sum_{k=0}^{K} C_A \cdot \left(
  \big(2\,C_I\,G_K(\mu)\big)^{\frac{1}{\mu}} \,
\big(F_k(s)\big)^{-\frac{1}{\mu+1}}\,\epsilon^{-\frac{1}{\mu}} +1 \right)
F_k(s) \, \TolXkm1\\
&\le C\cdot \big(G_K(\mu)\big)^{\frac{1}{\mu}}\epsilon^{-\frac{1}{\mu}}
\sum_{k=0}^{K} \big(F_k(s)\big)^{\frac{\mu}{\mu+1}}\,\TolXkm1
+ C_A\, \sum_{k=0}^{K}F_k(s)\,\TolXkm1\\
&= C\cdot \big(G_K(\mu)\big)^{\frac{\mu+1}{\mu}} \, \epsilon^{-\frac{1}{\mu}}
+ C_A\, \sum_{k=0}^{K}F_k(s)\,\TolXkm1
\end{align*}
with $C=C_A\,\big(2\,C_I\big)^{\frac{1}{\mu}}$.

The optimal tolerances $\TolYk$ can be directly determined from
assumption (A2). Thus with~$M_k$ defined above we get
\begin{equation}
\TolYKmk = \big(2\,C_Y\,G_K(\mu)\big)^{-1}\,\big(F_k(s) \big)^{\frac{\mu}{\mu+1}}\,
\TolXkm1\,\epsilon.
\end{equation}
Note that with these values $\sum_{k=0}^{K} C_Y\cdot \TolYk=\epsilon/2$,
which gives the desired accuracy in (\ref{est-stocherr}).\qedwhite

Observe that the function $G_K(\mu)$ as well as the second term in
\eqref{epscost} still depend on $\epsilon$, because
$K$ is a function of $\epsilon$. In this way, the choice of the
tolerances $\TolXk$ has an influence on the rate $-1/\mu$, which could
be further optimized. However, in the following we will  restrict our attention to a typical geometric
design with $\TolXk=q^k\,\TolXzero,\;k=1,2,\ldots,$ with a positive reduction
factor $q<1$.
The overall cost can then be estimated using a standard construction,
leading to the following result.
\begin{theorem}
\label{th:c_eps}
Let the sequence of spatial tolerances $\{\TolXk\}_{k=0,1,\ldots,K}$
in (\ref{seq-tolx}) be defined by $\TolXk=q^k\,\TolXzero$ with a reduction
factor $q<1$. Suppose assumptions (A1) and (A2) hold. Then, for any
$\epsilon<1$, there exists an integer $K=K(\epsilon)$ such that
\begin{equation}
\| u - u_K^{(\ML)}\|_{L^2_\rho(\Gamma;H^1_0(D))} \le \epsilon
\end{equation}
and
\begin{equation}
\ds C_\epsilon^{(\ML)} \lesssim \left\{
\begin{array}{ll}
\epsilon^{-\frac{1}{\mu}}                                 & \mbox{if } s\mu<1\\[1mm]
\epsilon^{-\frac{1}{\mu}}|\log\epsilon |^{1+\frac{1}{\mu}}& \mbox{if } s\mu=1\\[1mm]
\epsilon^{-s}                                             & \mbox{if } s\mu>1.
\end{array}
\right.
\end{equation}
\end{theorem}
\noindent {\it Proof:}
See \cite{Teckentrup2013,TeckentrupJantschWebsterGunzburger2012}. \qedwhite\pagebreak

In typical applications, it is usually more natural to consider a functional $\psi$
of the solution $u$ instead of the solution itself. Thus, suppose
a (possibly nonlinear) functional $\psi:H^1_0(D)\rightarrow\R$ with
$\psi(0)=0$ is given. In this case, we define
the following single-level and multi-level stochastic collocation approximations:
\begin{eqnarray}
\psi_{K}^{(\SL)} &\!:=\!& \cI_{M_K} \left[ \psi(u_K)\right], \\
\psi_{K}^{(\ML)} &\!:=\!& \sum_{k=0}^{K} \cI_{M_{K-k}}
\left[ \psi(u_k) - \psi(u_{k-1})\right]
\end{eqnarray}
with $u_{-1}:=0$. As in (\ref{equ-space-err}) and (\ref{equ-stoerr}),
we can ensure that for the adaptive error control of the expected
values, for all $k=0,\ldots,K,$ we have
\begin{eqnarray}
\label{qoi_cx}
\left| \E \left[ \psi(u) - \psi(u_k) \right]\right| &\!\le\!& C_X\cdot \TolXk\,,\\[2mm]
\label{qoi_cs}
\left| \E \left[ \psi(u_k) - \psi(u_{k-1})
-\cI_{M_{K-k}} \left[ \psi(u_k) - \psi(u_{k-1})\right]\right]\right| &\!\le\!&
C_Y\cdot \TolYKmk \,.
\end{eqnarray}
In practice, the following analogue of  Theorem \ref{th:c_eps} holds
for the expected value of the error of the multilevel approximation of functionals.
\begin{proposition}
\label{th:qoi}
Let the sequence of spatial tolerances $\{\TolXk\}_{k=0,1,\ldots,K}$
in (\ref{seq-tolx}) be defined by $\TolXk=q^k\,\TolXzero$ with a reduction
factor $q<1$. Suppose assumptions (A1) and (A2) hold with convergence
rates $\mu^*$ and $s^*$.
Then, for any
$\epsilon<1$, there exists an integer $K(\epsilon)$ such that
\begin{equation}
\left| \E \left[ \psi(u) - \psi_K^{(\ML)} \right] \right| \le \epsilon
\end{equation}
and
\begin{equation}
\ds C_\epsilon^{(\ML)} \lesssim \left\{
\begin{array}{ll}
\epsilon^{-\frac{1}{\mu^*}}                                 & \mbox{if } {s^*}{\mu^*}<1\\[2mm]
\epsilon^{-\frac{1}{\mu^*}}|\log\epsilon |^{1+\frac{1}{\mu^*}}& \mbox{if } {s^*}{\mu^*}=1\\[2mm]
\epsilon^{-{s^*}}                                             & \mbox{if } {s^*}{\mu^*}>1.
\end{array}
\right.
\end{equation}
\end{proposition}
\noindent {\it Proof:} see the discussion of
 Proposition 4.6 in~\cite{TeckentrupJantschWebsterGunzburger2012}.

Note that in analogy to strong and weak convergence of solutions
of stochastic differential equations, one would expect the error in mean
associated with approximating a functional $\psi(u)$ to decrease at
a faster rate than the error in norm associated with approximating
the entire solution $u$. Thus, we anticipate that $\mu^*>\mu$.
Moreover, the convergence rate of the error
in the functional $\psi(u_k)$ is in general larger than the
convergence rate of the error in the $H_1$-norm, which leads to a
smaller value $s^*<s$ in assumption (A1).
Specifically, in the case of $H^2$-regularity in space and using
an optimal linear solver, such as multigrid, we anticipate that ${s^*}=d/2$.

We conclude this section with a complete algorithmic description of our adaptive
multilevel stochastic collocation method. Table~\ref{tab-algo-amlsc} illustrates the main steps.
The strategy is  self-adaptive in nature. Thus, once the tolerances $\{\TolXk\}_{k=0,\ldots,K}$ and
$\{\TolYk\}_{k=0,\ldots,K}$ are set, the algorithm delivers
an approximate functional $\psi_K^{(\ML)}$ with
accuracy close to the user-prescribed tolerance~$\epsilon$. An
analogous algorithm can be defined to deliver the numerical
solution $u_K^{(\ML)}$ close to a user-prescribed tolerance
$\epsilon$.
A priori, to obtain optimal results, the reliability of the estimation
for the adaptive spatial discretization and the adaptive Smolyak algorithm needs
to be studied in order to provide values for the constants $C_X$ and $C_Y$.
More specifically, while the spatial tolerances $\TolXk$ can be freely chosen
(by fixing the number of levels $K$
and the reduction factor  $q$), the optimal choice of the tolerances
$\TolYk$ in (\ref{opt-toly}) requires estimates of the parameters
$s^*$ and $\mu^*$. A discussion of how to effect this in a preprocessing step
is given in the next section. A crucial point to note is that, even without this information, the adaptive anisotropic
Smolyak algorithm will automatically detect the importance of various directions in the parameter space
$\Gamma\subset\R^N$.

\begin{table}[t!]
\centering
\begin{tabular}[t]{ll}
\hline
\multicolumn{2}{l}{Algorithm: Adaptive Multilevel
Stochastic Collocation Method} \\
\hline\rule{0mm}{5mm}\hspace{-0.1cm}
1. & Given $\epsilon$ and $q$, \;estimate $C_X$, $C_Y$,
     $s^*$, $\mu^*$ and $K=K(\epsilon)$.\\[0.5ex]
2. & Set coarsest spatial tolerance: \  $\TolXzero:=\epsilon/\big(2C_X\,q^K\big)$.\\[0.5ex]
3. & Set other spatial tolerances: \ $\TolXk:=q^k\,\TolXzero$, $k=1,\ldots,K$.\\[0.5ex]
4. & Compute $\TolXm1:=\E[\cI_{M_0}[\psi(u_0)]]$ with
     tolerances $\TolXzero$ and $\TolYzero :=\TolXzero$.\\[0.5ex]
5. & Set stochastic tolerances:\\
   & ~~$\TolYKmk := \big(2C_YG_K(\mu^*)\big)^{-1} \,
     \big(F_k(s^*)\big)^{\mu^*/(\mu^*+1)} \, \TolXkm1
     \, \epsilon$,
     $k=0,\ldots,K$, \\[0.8ex]
    & ~~with $F_k$ and $G_K$ defined in (\ref{Fkdef}) and
      (\ref{Gkdef}), respectively.\\[0.5ex]
6. & Coarsest level $k=0$ (reusing samples from Step 4):\\[0.2ex]
   & ~~Compute $E_0:=\E[\cI_{M_K}[\psi(u_0)]]$ with tolerances $\TolXzero$ and $\TolYK$.\\[0.5ex]
7. & Multilevel differences $k=1,\ldots,K$ (reusing samples from level $k-1$):\\[0.2ex]
   & ~~Compute $E_k:=\E[\cI_{M_{K-k}}[\psi(u_{k})-\psi(u_{k-1})]]$
     with tolerances
     $\TolXk$, $\TolXkm1$ and $\TolYKmk$.\\[0.5ex]
8. & Compute $\E\big[\psi^{(\ML)}_K\big]:=\sum_{k=0,\ldots,K}E_k$.\\[1mm]
\hline
\end{tabular}\\
\parbox{13cm}{
\caption{\small Algorithm to approximate solution functionals $\psi(u)$ by an
adaptive multilevel stochastic collocation method.}
\label{tab-algo-amlsc}
}
\end{table}

An important point already mentioned in
\cite{TeckentrupJantschWebsterGunzburger2012} is that the optimal
rounded values for the number of samples, $M_k$, will not be used
by the algorithm, because they do not necessarily correspond to an
adaptive sparse grid level. However, for each level $k$, the tolerance
$\TolYk$ can be ensured by choosing $\widetilde{M}_k\ge M_k$ slightly larger,
resulting in a slight inefficiency of the sparse grid approximation.
Note that, in practice, the same behaviour is observed for adaptive
spatial discretizations. In any case, no restart is needed; cf.~\cite[Section 6.3]{TeckentrupJantschWebsterGunzburger2012}.

\section{Numerical examples}
\label{sec:numerics}

First, we provide general information on the adaptive components used. Then,
numerical results are presented for two isotropic diffusion problems with uncertain source
term and uncertain geometry, respectively. All calculations have been done
with \matlab\ version R2017a on a Latitude $7280$ with an i5-7300U Intel
processor running at 2.7 GHz.

For the spatial approximation considered in the two examples,
we use the adaptive piecewise linear finite element method implemented by Funken, Praetorius and Wissgott in the \matlab\  package {\sl p1afem}.\footnote{The {\sl p1afem} software package
can be downloaded from the author's webpage~\cite{FunkenPraetoriusWissgott}.}
Using \matlab\ built-in functions and vectorization for an efficient realization, the code performs with almost linear complexity in terms of degrees of freedom with respect to the runtime.
A general description of the underlying ideas can be found in \cite{FunkenPraetoriusWissgott2011}
and the technical report \cite{FunkenPraetoriusWissgott2008} provides  detailed documentation. The code is easy to modify. In
order to control the accuracy of solution functionals,
the dual weighted residual method (DWRM) introduced by Becker \&
Rannacher~\cite{BeckerRannacher2001} is adopted. In what follows, we will give a short summary
of the underlying principles that are relevant for our implementation.

Let $y\in\Gamma$ be fixed. For ease of presentation we write $u=u(\cdot,y)$
and consider the variational problem (\ref{eqn_rand1}) with $a(u,v) = \int_D \nabla u \cdot \nabla v  \,dx$.
We also suppose the solution functional $\psi(u)$ takes the specific form
\begin{equation}
\psi(u) = \int_D u^2 \,dx ,
\end{equation}
as discussed in~\cite[Example 5.13]{TeckentrupJantschWebsterGunzburger2012}.
Let $u_k\in V_k$ be the finite element
solution computed on the (adaptive) mesh $\mathcal{T}_k$. Then the DWRM provides a
representation of the error in the solution functional in the form
\begin{align}  \label{dwrm}
 \psi(u) - \psi(u_k)
 &=  \int_D 2 u_k (u -u_k)  \,dx  + \int_D (u -u_k)^2  \,dx
  \approx \int_D 2 u_k (u -u_k)  \,dx ,
\end{align}
by simply neglecting the higher-order term.
Letting  $w$ be the exact solution of the linearised dual problem
\begin{equation} \label{dualprob}
\int_D \nabla v \cdot  \nabla w \,dx =  \int_D 2 u_k  v \,dx \quad \forall v\in H^1_0(D)
\end{equation}
and letting $w_k \in V_k$ be the finite element approximation of the dual solution on the same mesh, we have, using Galerkin orthogonality,
\begin{align}
 \psi(u) - \psi(u_k)
& \approx \int_D 2 u_k (u -u_k)  \,dx  =  \int_D \nabla (u -u_k) \cdot  \nabla (w -w_k) \,dx \\
&=  \sum_{T\in\mathcal{T}_k}
\left\{
\int_T f(x,y) (w-w_k) - \nabla u_k \cdot \nabla (w-w_k)\,dx \right\} .
\end{align}
In practice, when solving the dual problem we
compute an approximation $\phi_k\approx w-w_k$ of its error in the hierarchical surplus
space of quadratic finite elements. Hierarchical error estimators using this approach are implemented in {\sl p1afem} for the primal solution $u_k$.

Putting this all together, we obtain
\begin{equation}
\psi(u) - \psi(u_k) \approx \sum_{T\in\mathcal{T}_k} \eta_T
\quad\hbox{with}\quad
 \eta_T = \int_T f(x,y)\,  \phi_k - \nabla u_k \cdot \nabla \phi_k\,dx.
\end{equation}
We use $|\eta_T|$ as refinement indicators and
mark elements $T\in\mathcal{T}_k$ for refinement using the standard D\"orfler criterion
from~\cite{Doerfler1996}, which determines the minimal set $\mathcal{M}\subset\mathcal{T}_k$
such that $\theta\sum_{T\in\mathcal{T}_k}|\eta_T|\le \sum_{T\in\mathcal{M}}|\eta_T|$. We
typically set a value of $\theta=0.6$ in our calculations. Refinement
by newest vertex bisection is applied to guarantee nested finite element spaces and
the optimal convergence of the adaptive finite element method, see~\cite{KreuzerSiebert2011}.
The adaptive process is terminated when the absolute value of $\sum_{T\in\mathcal{T}_k} \eta_T$
is less than a prescribed tolerance.

Turning to the adaptive anisotropic Smolyak algorithm, the main idea for the construction of the sparse grid interpolation operators in (\ref{stoch-intpl-op}) is to use the hierarchical decomposition
\begin{equation}
\label{sparse-grid-intpl}
\cI_{M_k}[u_h](y) = \sum_{\bfi\in I} \triangle^{m(\bfi)}[u_h](y)
:=\sum_{\bfi\in I} \bigotimes_{n=1}^{N}
\left( \cI_n^{m(i_n)}[u_h](y) - \cI_n^{m(i_n-1)}[u_h](y)\right)
\end{equation}
with multi-indices $\bfi=(i_1,\ldots,i_N)\in I\subset\N^N_+$,
$m(\bfi)=(m(i_1),\ldots,m(i_N))$, and
univariate polynomial interpolation operators
$\cI_n^{m(i_n)}:C^0(\Gamma_n)\rightarrow\PO_{m(i_n)-1}$, which
use $m(i_n)$ collocation points to construct a polynomial interpolant
in $y_n\in\Gamma_n$ of degree at most $m(i_n)-1$. The operators
$\triangle^{m(\bfi)}$ are often referred to as hierarchical
surplus operators. The function
$m$ has to satisfy $m(0)=0$, $m(1)=1$, and $m(i)<m(i+1)$. We set
$\cI^0_n=0$ for all $n=1,\ldots,N$ and
use the nested sequence of univariate Clenshaw--Curtis nodes with
$m(i)=2^{i-1}+1$ if $i>1$. In (\ref{sparse-grid-intpl}), $M_k$ is
then the number of all explored quadrature points in $\Gamma$
determined by $m(\bfi)$. To get
good approximation properties, the index set $I$ should satisfy the
downward closed set property, i.e.,
\begin{equation}
\label{downward-closed}
\text{if }\bfi\in I,\text{ then }\bfi-\bfe_j\in I \text{ for all }
j=1,\ldots,N \text{ such that }i_j>1.
\end{equation}
As usual, we ensure that ${\bf 1}\in I$ to also recover
constant functions.

The hierarchical structure in (\ref{sparse-grid-intpl}) allows to
interpret updates that are derived by adding further differences
$\triangle^{m(\bfi_a)}$, i.e., enhancing the index set $I$ by an
admissible index $\bfi_a$ that satisfies \eqref{downward-closed},
as error indicators for already
computed approximations. There are several adaptive strategies
available. One could explore the whole margin of $I$ defined by
\begin{equation}
M_I := \{ \bfi\in\N^N_+\backslash I:\,\bfi-\bfe_n\in I\text{ for some }
n\in\{1,\ldots,N\}\}.
\end{equation}
Generally, this approach is computationally challenging and yields
a fast increase of quadrature points. Instead, as suggested by
Gerstner \& Griebel~\cite{GerstnerGriebel2003}, the margin
is reduced to the set
 \begin{equation}
R_I := \{ \bfi\in M_I:\,\bfi-\bfe_n\in I\text{ for all }
n=1,\ldots,N \text{ with }i_n>1\}.
\end{equation}
In each step, the adaptive Smolyak algorithm
computes the profits $\triangle^{m(\bfi_a)}$ for all
$\bfi_a\in R_I$ -- reusing already computed profits --
and replaces the index in $R_I$ with the highest profit,
say $\bfi_{max}$, by its admissible neighbours taken from the set
$\{\bfi_{max}+\bfe_j, j=1,\ldots,N\}$. These neighbours are then
explored next. The algorithm stops if the absolute value of the
highest profit is less than a prescribed tolerance. This
adaptive strategy is implemented in \matlab\ in the
{\sl Sparse Grid Kit}\footnote{This package can be downloaded from
the CSQI website~\cite{TamelliniNobileGuignardTeseiSprungk2017}.}
and its numerical performance is discussed
in the review paper~\cite{BeckNobileTamelliniTempone2011}.

We adopt this refinement strategy with
one minor deviation in our numerical experiments: namely, at the final iteration step, instead of adding the new profits to the interpolant, the previous approximate value is returned. In that way, the last highest profit can  be
considered as a more realistic error indicator. This adaptation
allows a better understanding of the convergence behaviour of
the multilevel approach. In practical calculations one
would almost certainly use the final value.

\subsection{Uncertain source and boundary conditions: $d=2$, $N=2$}
We will refer to this example as the {\it one peak\/} test problem. It was
introduced by Kornhuber \& Youett~\cite[Sect.~5.1]{KornhuberYouett2018}
in order to  assess  the efficiency of  adaptive multilevel Monte Carlo methods.
The challenge is  to solve the Poisson equation $-\nabla^2 u = f $ in a unit square domain
$D=(-1,1)\times(-1,1)$ with Dirichlet boundary data  $u=g$ on $\partial D$ (points in $D$ are
 represented by ${\vfld x}=(x_1,x_2)$).  The source term $f$ and boundary data are  {\it uncertain} and  are
parameterised by  $ \bfy =(y_1,y_2)$, representing the image  of  a pair
of independent random  variables  with  $y_j \sim {\cal U}[-1/4,1/4]$.
In the isotropic case studied in~\cite{KornhuberYouett2018},
the source term $f$ and boundary data $g$ are chosen so that the uncertain PDE problem
has a specific pathwise solution given by
$$
u({\vfld x}, \bfy) = \exp ( - \beta  \{ (x_1 -y_1)^2 + (x_2 -y_2)^2 \} ),
$$
where the scaling factor $\beta>0$ is  chosen to be large---so as to generate
a highly localised Gaussian profile centered at the  uncertain spatial location $(y_1,y_2)$.
The value of $\beta$ that we will take  in this study is $\beta=50$. (The other values discussed
in~ \cite{KornhuberYouett2018} are $\beta=10$ and $\beta=150$.)

In this work, the test problem  is made {\it anisotropic} by  scaling the solution in  the first coordinate
direction  by a linear function $\alpha(y_1)= 18 y_1 + 11/2$ so that
the $\alpha$ takes values in the interval $[-1,10]$. The corresponding pathwise solution is then given by
$$
u({\vfld x},\bfy) = \exp ( - 50  \{  \alpha(y_1) (x_1 -y_1)^2 + (x_2 -y_2)^2 \} ).
$$
The goal is to approximate the following quantity of interest (QoI)
\begin{equation}
\E [\psi(u)] = \E \left[ \int_D u^2(x,y)\,dx\right].
\end{equation}

\begin{figure}[t]
\centering
\includegraphics[width=0.9\textwidth]{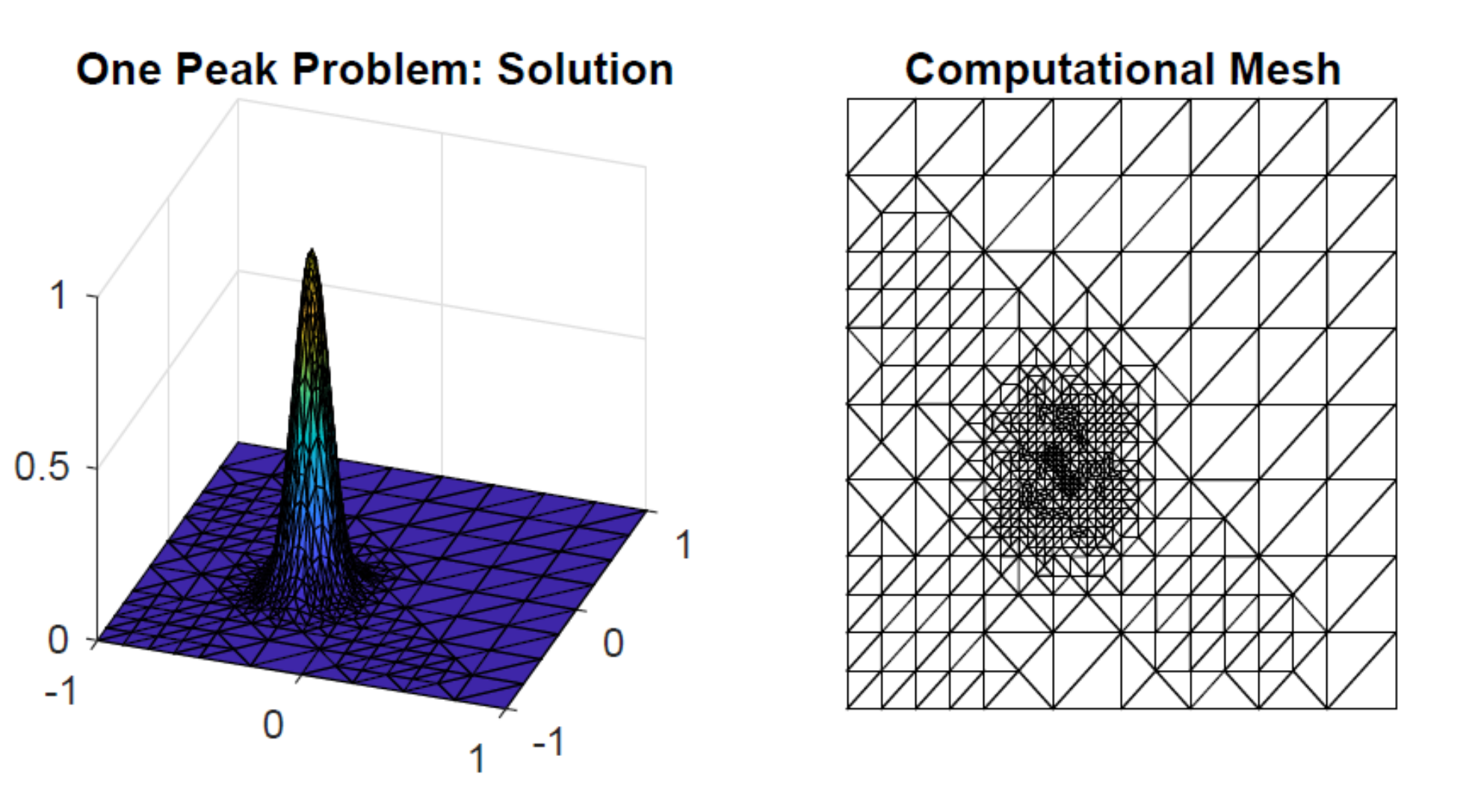}
\parbox{13cm}{
\caption{\small One peak problem: solution for $y=(-0.22,-0.22)^T$ (left)
and corresponding adaptive mesh generated using the
dual weighted residual approach
and a tolerance of $\TolX=2\times 10^{-3}$ (right).}}
\label{fig:OnePeakSolMesh}
\end{figure}

In our experiments we solve the simplified  problem with $u=0$ on $\partial D$. (This has no
impact on accuracy for the error tolerances that we will considered; see the discussion in  the appendix).
The initial mesh with $81$ vertices is generated by taking three uniform refinement steps, starting
from two triangles with common edge from $(-1,-1)$ to $(1,1)$. For the specific
parameter choice $y=(-0.22,-0.22)^T$ and spatial tolerance $\TolX=2\times10^{-3}$,
the adaptive algorithm terminated after 5 steps giving
the numerical solution and  corresponding mesh shown in Fig.~\ref{fig:OnePeakSolMesh}.

\begin{figure}[t!]
\centering
\hspace*{-0.75cm}\includegraphics[width=0.5\textwidth]{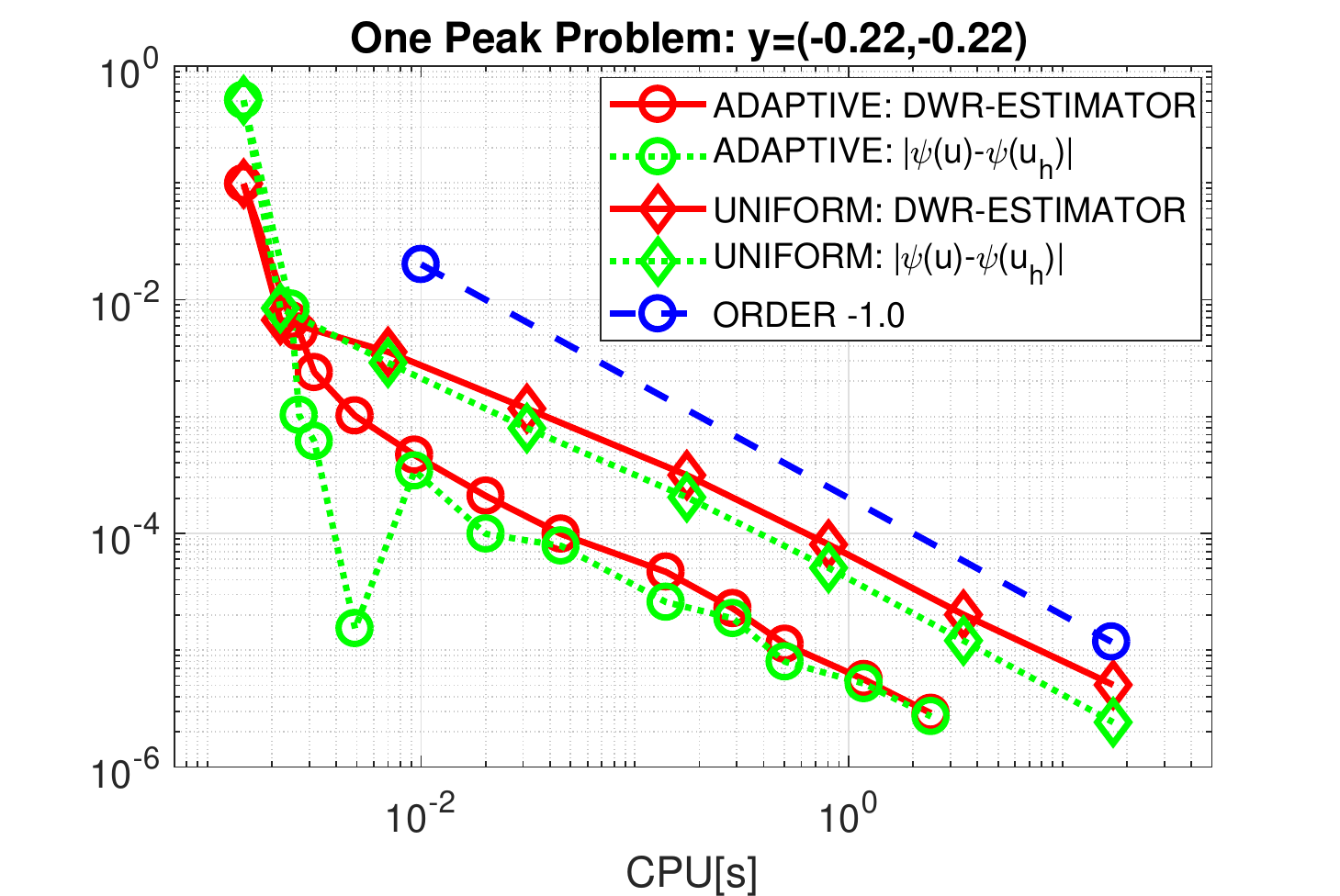}
\includegraphics[width=0.5\textwidth]{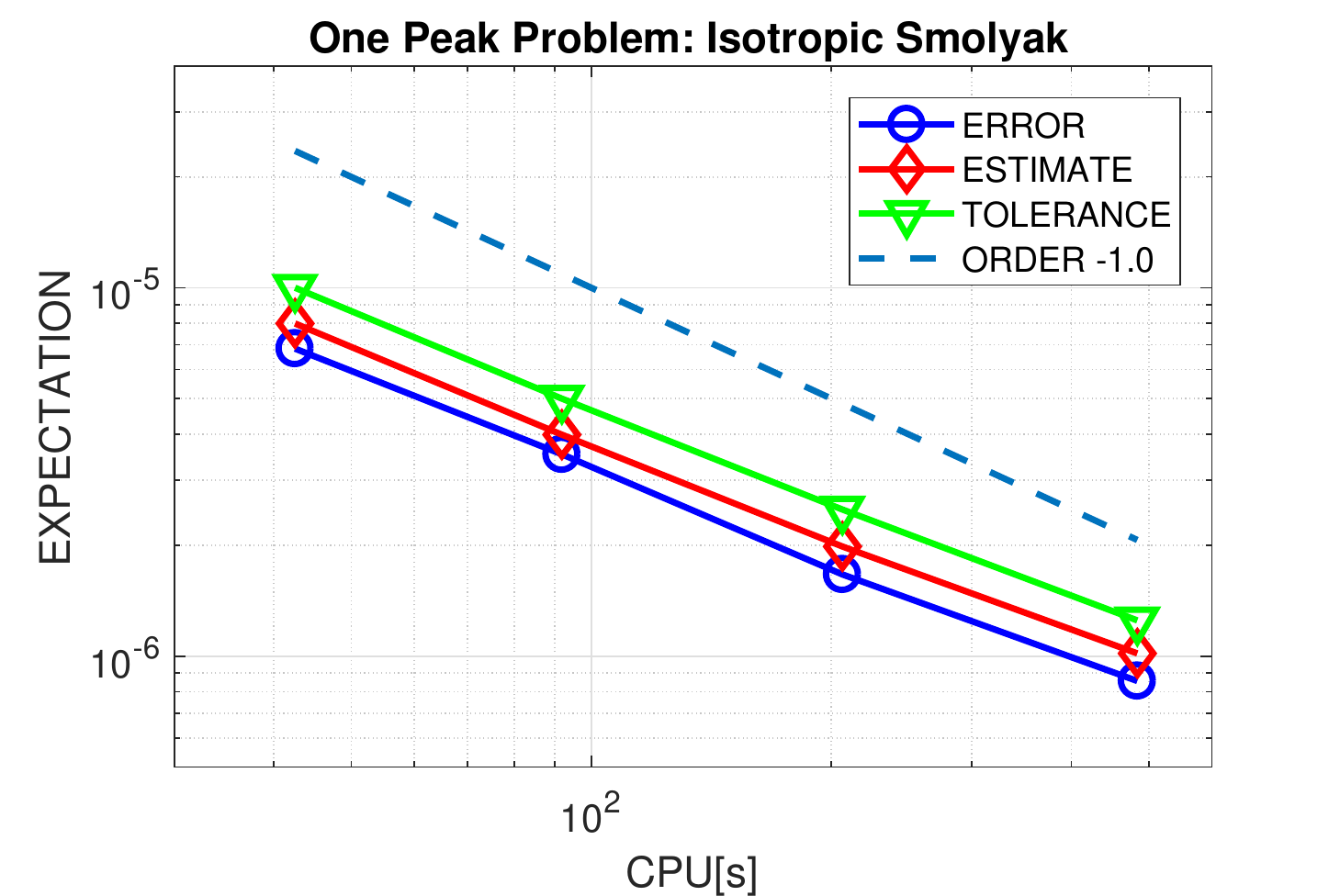}
\parbox{13cm}{
\caption{\small One peak problem: history of error estimators and exact errors obtained by adaptive and uniform spatial refinements for
$y=(-0.22,-0.22)^T$
(left); history of error estimators and
exact errors obtained by adaptive spatial refinements averaged
over the parameter domain using an isotropic Smolyak approximation
with $145$ collocation points and spatial tolerances $\TolX=10^{-5}/2^i$, $i=0,\ldots,3$ (right). The numerically observed convergence
orders for $\psi(u)$ and $\E[\psi(u)]$ in terms of CPU time
are in all cases close to $-1$.}
\label{fig:OnePeakSolEst}
}
\end{figure}

Running the adaptive algorithm with a tighter tolerance of $\TolX=5\times10^{-6}$
generated a concentrated mesh with $158\,734$ points after 13 refinement steps.
The left part of  Fig.~\ref{fig:OnePeakSolEst} shows the efficiency of the locally
adaptive procedure in comparison with a uniform refinement strategy.
For smaller tolerances, the gain of efficiency in terms of overall cpu time
is a factor $10$. Comparing with a reference value $\psi(u)= 0.025315675 \ldots$,
generated by running the adaptive algorithm with a very small tolerance, we see that
the effectivity index, i.e., the ratio between estimator and true error,
tends asymptotically to $1$ for both uniform and adaptive refinement.

The very good quality of the spatial error estimation process
is still maintained if we sample over the whole parameter space
using an isotropic sparse grid of $145$ collocation points
-- generated using {\sl Sparse Grid Kit} by doubling of points in each stochastic
dimension. In this case we have a reference value
of $\E [\psi(u)] = 0.015095545 \ldots$, generated analytically using an
asymptotic argument---details are given in the appendix.
The error estimates for adaptive meshes associated  with
different spatial  error tolerances are plotted in
the right part of Fig.~\ref{fig:OnePeakSolEst}.
Observe that the estimators deliver upper bounds for the numerical errors
and the tolerances are always satisfied.
As expected from the theory, the convergence rates for $\psi(u)$ and
$\E[\psi(u)]$ in terms of computing time are both close to $-1$. So we have
$C_X\!=\!1$ in (\ref{qoi_cx}) and $s^*\!=\!1$ in assumption (A1)
in Section~\ref{methodology}.

\begin{figure}[h!]
\centering
\includegraphics[width=0.48\textwidth]{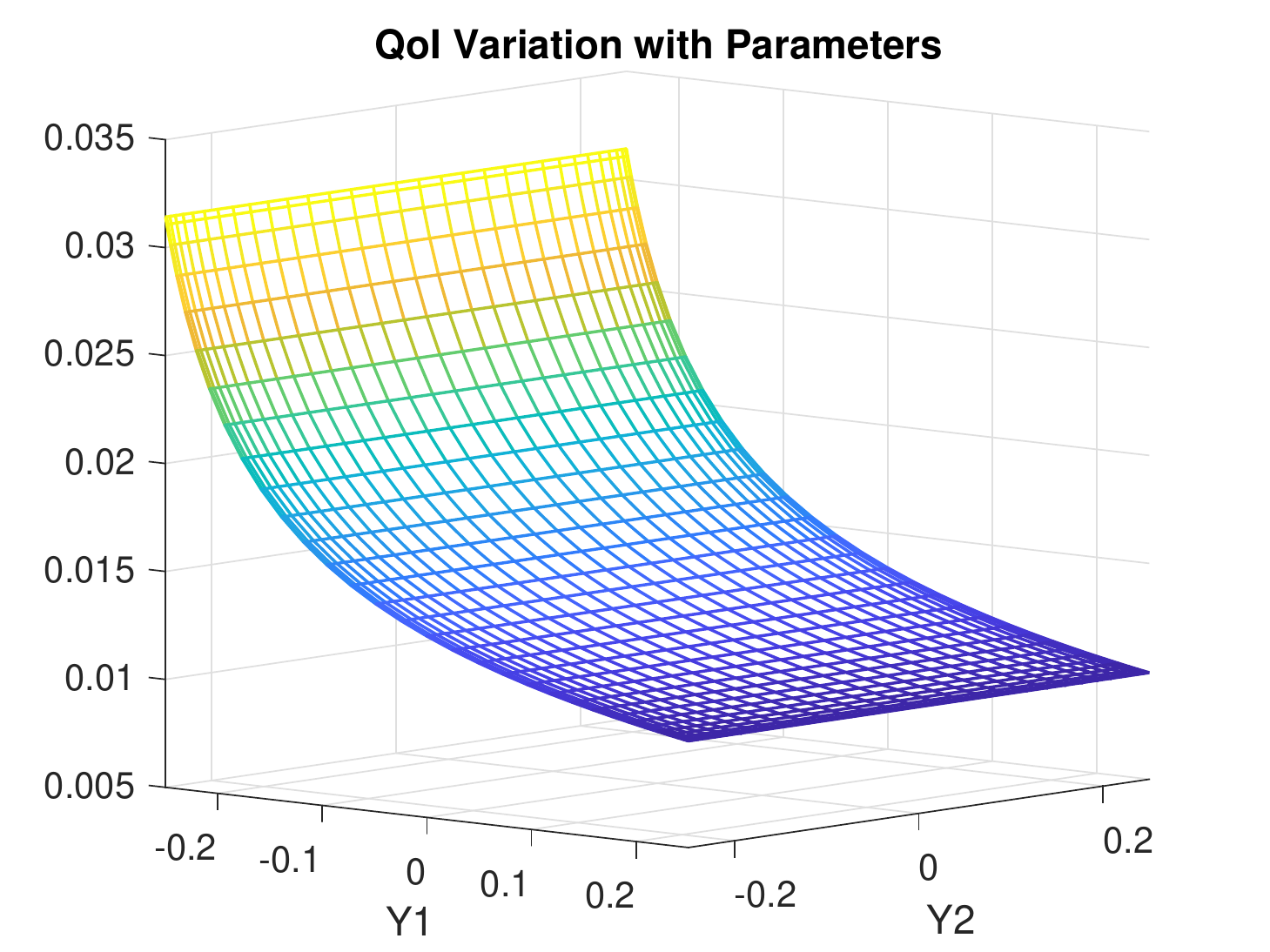}
\includegraphics[width=0.5\textwidth]{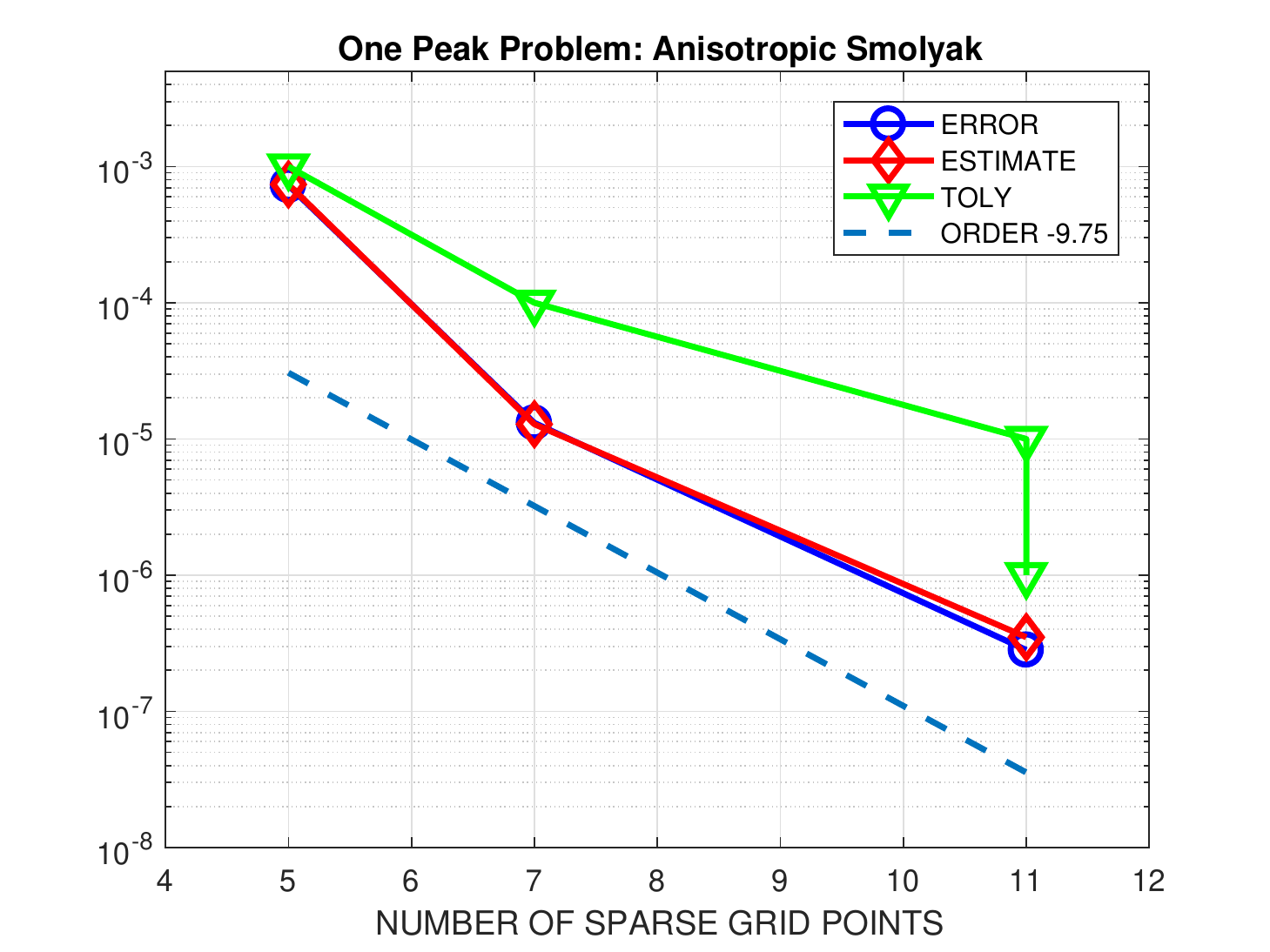}
\parbox{13cm}{
\caption{\small One peak problem: reference solution computed
using  $33\times 33$ collocation points (left);
convergence of error estimators and exact errors
for anisotropic Smolyak approximations with stochastic tolerances
$\TolY=10^{-(i+3)},i=0,\ldots,3$ (right).
The estimated convergence
order for $\E[\psi(u)]$ in terms of collocation points is $-9.75$.}
\label{fig:OnePeakCollEst}
}
\end{figure}

Next we consider the convergence behaviour and the quality of the error
estimates for the adaptive {\it anisotropic} Smolyak algorithm.
A reference solution generated by solving the problem with a small spatial error tolerance on a $33\times 33$ tensor-product grid of Clenshaw--Curtis points is shown in the left part of Fig.~\ref{fig:OnePeakCollEst}. Looking at the solution structure, it is evident that the QoI has inherent structure  that can be exploited by an adaptive
strategy: more specifically, very few sample points are needed to accurately compute it$\,$!
The right part of Fig.~\ref{fig:OnePeakCollEst}
shows the results for a decreasing sequence of
stochastic tolerances $\TolY=10^{-3},\ldots,10^{-6}$, where we have fixed the spatial tolerance $\TolX=10^{-7}$. The corresponding
numbers of collocation points are $5,7,11$, and $11$. Looking
at these more closely we find that all the quadrature rules generated
by the adaptive procedure are one-dimensional rules in the $y_1$ direction matching the variation in the reference.
Other features of note are that the prescribed tolerances
are always satisfied and that the errors (with respect to the reference
solution) are nicely estimated by the error indicators.
A least-squares fit gives an averaged
value of $\mu^*=9.75$ for the convergence order
so that assumption (A2) in Section~\ref{methodology} is
satisfied. An approximate value $\mu^*\approx 10.74$ can
be computed with a much lower spatial tolerance $\eta_X=1.25\times 10^{-4}$ and
$19$ collocation points in a few seconds.

Having estimated the parameters in the second step of the algorithm
in Table~\ref{tab-algo-amlsc},
we now run the adaptive multilevel approach with overall accuracy
requirements of $\epsilon=10^{-5},5\times10^{-6},2.5\times10^{-6}$ and
$10^{-6}$.
To illustrate the performance gains, we simply set
$K\!=\!2$ (so that three levels are used) and assign a spatial error   reduction factor of
$q=0.2$. The spatial tolerances are then given by
$\TolXk=\epsilon q^{k-2}/2$ with $k=0,1,2$. To calculate, in a first
step, a sufficiently accurate approximation for the tolerance
$\TolXm1=\E[\psi(u_0)]$ at reasonable cost, we apply
the anisotropic Smolyak algorithm with $\TolX=\TolY=\TolXzero$.
Note that these samples can be reused later in the first level of the
multilevel scheme. Eventually,
the stochastic tolerances are derived from (\ref{opt-toly}) with
$C_Y=0.1$, $\mu^*=9.75$ and $s^*=1$.

 \begin{figure}[t]
\centering
\includegraphics[width=0.8\textwidth]{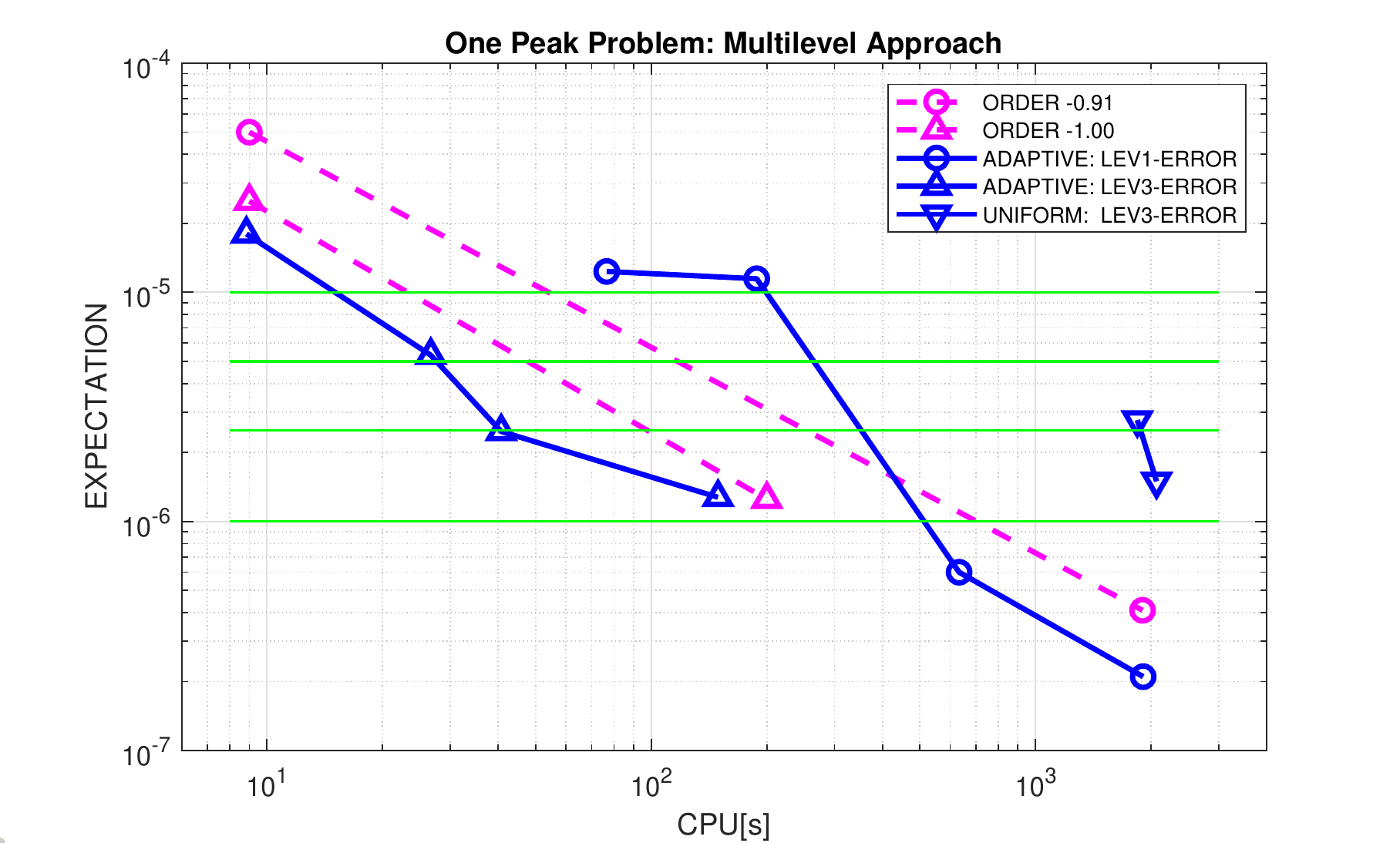}
\parbox{13cm}{
\caption{\small One peak problem: errors for the expected values $\E[\psi_{2}^{(\ML)}]$
and $\E[\psi_{2}^{(\SL)}]$ for the three-level (blue triangles) and the one-level
(blue circles) approach with adaptive spatial meshes for overall
accuracy requirements
$\epsilon=10^{-5},5\times10^{-6},2.5\times10^{-6},10^{-6}$, respectively (green lines). The orders of
convergence predicted by Theorem~\ref{th:qoi} are $-1$ and $-0.91$,
respectively (dashed magenta lines).}
\label{fig:OnePeakMultiLevel}
}
\end{figure}

Results for the three-level and single-level approach are summarised in
Fig.~\ref{fig:OnePeakMultiLevel}. For adaptive spatial meshes,
we observe that the errors of the
expected values $\E[\psi_{2}^{(\ML)}]$ are very close to the
prescribed tolerances. This is not always the case for the values
$\E[\psi_{2}^{(\SL)}]$ computed by the single-level approach.
The three-level approach performs reliably and
outperforms the single-level version clearly. The
orders of convergence, $p_{\ML}\!=\!-1/s^*\!=\!-1$ and
$p_{\SL}\!=\!-1/(s^*+1/{\mu^*})\!\approx\!-0.91$,
predicted by Theorem \ref{th:qoi} for the accuracy in terms of computational
complexity are also visible. Already for coarse tolerances, the three-level approach
with uniform spatial meshes
shown plotted with inverted blue triangles) is not competitive.
For higher tolerances, the calculation exceeded
memory requirements (more than $1.6\times10^{7}$ FE degrees of freedom).

These CPU timings are impressive.
We have also solved the one peak test problem using  an efficient adaptive stochastic Galerkin (SG) approximation strategy.  While the linear algebra
associated with the Galerkin formulation is decoupled in this case,
the computational overhead  of evaluating the right-hand-side vector
is a big limiting factor in terms of the relative efficiency. Using SG  the quantity
$$
I =
  \int_{\triangle}  \int_\Gamma  d({\vfld x},\bfy) \cdot f({\vfld x},\bfy) \,  \psi_{\nu_k}(\bfy)
\,  \phi_j({\vfld x}) \, d\bfy \, d\vfld x
$$
where
$$
f({\vfld x},\bfy) = d(x_1,x_2,y_1,y_2) \cdot  \exp ( - 50  \{  \alpha(y_1) (x_1 -y_1)^2 + (x_2 -y_2)^2 \} ),
$$
and
$$
d(x_1,x_2,y_1,y_2) = -10000 \left \{ \alpha^2(y_1)  (x_1 -y_1)^2 + (x_2 -y_2)^2 \right \}  + 100 (\alpha(y_1)+1)
$$
must be computed in every element $\triangle$ in the current subdivision
and for every parametric function $ \psi_{\nu_k}(\bfy)$ in the
active index set. This  is an extremely  demanding quadrature problem!

We have also tested an adaptive multilevel Monte Carlo method
(see, \cite{CliffeGilesScheichlTeckentrup2011,KornhuberYouett2018}) on this problem.
For the lowest tolerance, $\epsilon=10^{-5}$, and an average over $5$ independent
realizations, the three-level algorithm achieves an accuracy of $2.77\times10^{-6}$
in $7.76\times10^{4}$ sec. The numbers of averaged samples for each level
are $M_0=519634$, $M_1=6153$, and $M_2=243$. Obviously, the slow Monte
Carlo convergence
rate of $\mu=0.5$ is prohibitive for higher tolerances here.

\subsection{Uncertain geometry: $d=2$, $N=16$}
In our second example, we again consider a two-dimensional Poisson problem, but now with geometry of the computational domain being uncertain.
We will refer to this example as the {\it two hole\/} test problem.
The uncertain domain is defined by
\begin{equation}
D(y)=(0,6)\times (0,6)\backslash (D_1(y)\cup D_2(y)),
\end{equation}
where the holes $D_1(y)=\overline{P_1P_2P_3P_4}$ and
$D_2(y)=\overline{P_5P_6P_7P_8}$ are taken as quadrilaterals with
the uncertain vertices
\begin{equation}
\begin{array}{ll}
P_1=(1+y_1/a_1,1+y_2/a_2),&
P_2=(2+y_3/a_3,1+y_4/a_4),\\[1mm]
P_3=(2+y_5/a_5,3+y_6/a_6),&
P_4=(1+y_7/a_7,3+y_8/a_8),\\[1mm]
P_5=(4+y_9/a_9,1+y_{10}/a_{10}),&
P_6=(5+y_{11}/a_{11},1+y_{12}/a_{12}),\\[1mm]
P_7=(5+y_{13}/a_{13},5+y_{14}/a_{14}),&
P_8=(4+y_{15}/a_{15},5+y_{16}/a_{16}).
\end{array}
\end{equation}
Here, the random vector $y=(y_1,\ldots,y_{16})^T$ consists
of sixteen uniformly distributed random variables
$y_i\sim {\cal U} [-1,1]$, $i=1,\ldots,16$. We set
\begin{equation}
(a_1,\ldots,a_{16})=(5,5,5,5,5,5,5,5,10,10,10,10,20,20,20,20)
\end{equation}
to represent different strength of uncertainty and hence anisotropy
in the stochastic space. We impose a constant volume force $f\!\equiv\!1$ and
fix the component by homogeneous Dirichlet
boundary conditions on the whole boundary $\partial D(y)$, including
$\partial D_1(y)\cup\partial D_2(y)$ with stochastically
varying positions in space. Our goal is then to study the effect of
this uncertainty on the expectation of the overall displacement
calculated by
\begin{equation}
\E [\psi(u,y)] = \E \left[ \int_{D(y)} u^2(x,y)\,dx\right].
\end{equation}
This allows an assessment of the desired averaged load capacity of the component,
which takes into account uncertainties in the manufacturing process.

\begin{figure}[t!]
\centering
\hspace*{-1.35cm}\includegraphics[width=1.2\textwidth]{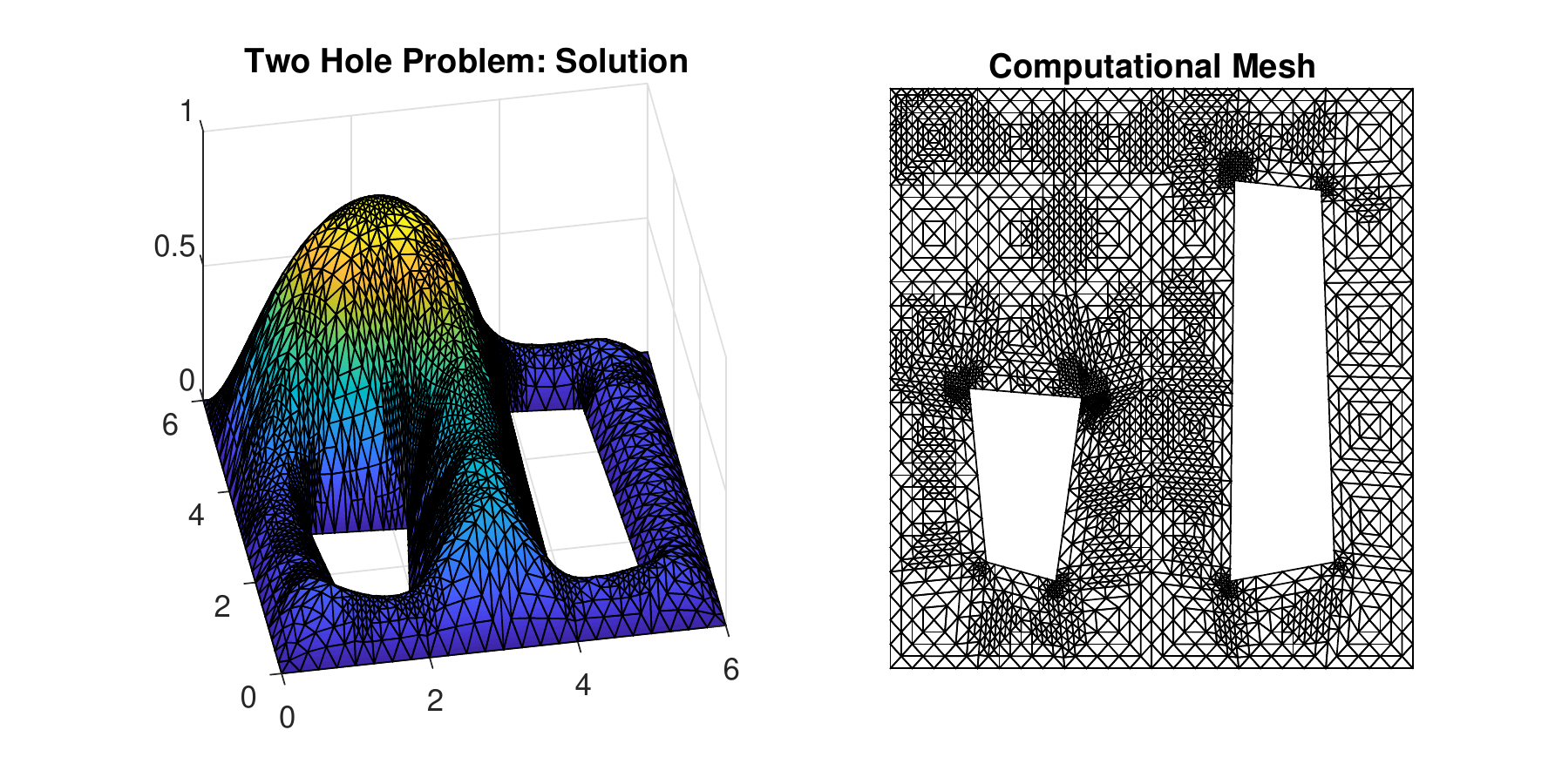}
\parbox{13cm}{
\caption{\small Two hole problem: Solution for the particular random
domain $D(y)$ with
$y=[0.5,0.5,-0.5,-0.5,-0.5,-0.5,1,-1,-1,-1,1,1,-1,1,-1,-1]$ (left)
and corresponding adaptive mesh based on the dual weighted residual
approach (right), using a spatial tolerance $\TolX=10^{-1}$ and
$2\,650$ mesh points to resolve the corner singularities.}
\label{fig:TwoHoleSolMesh}
}
\end{figure}

Applying a parameter-dependent map, each domain $D(y)$ can be mapped to the fixed
nominal domain $D_0=D(0)$ with $0\in\R^{16}$. Such a domain mapping approach was
introduced by Xiu \& Tartakovsky~\cite{XiuTartakovsky2007} and
allows us to reformulate the problem in the form (\ref{eqn_rand1}) with parameter-dependent
coefficients on $D_0$. The well established theory for elliptic partial differential
equations with random input data can then be applied without modifications to show the
well-posedness of the setting with random domains.

\begin{figure}[t]
\centering
\hspace*{-0.75cm}\includegraphics[width=0.55\textwidth]{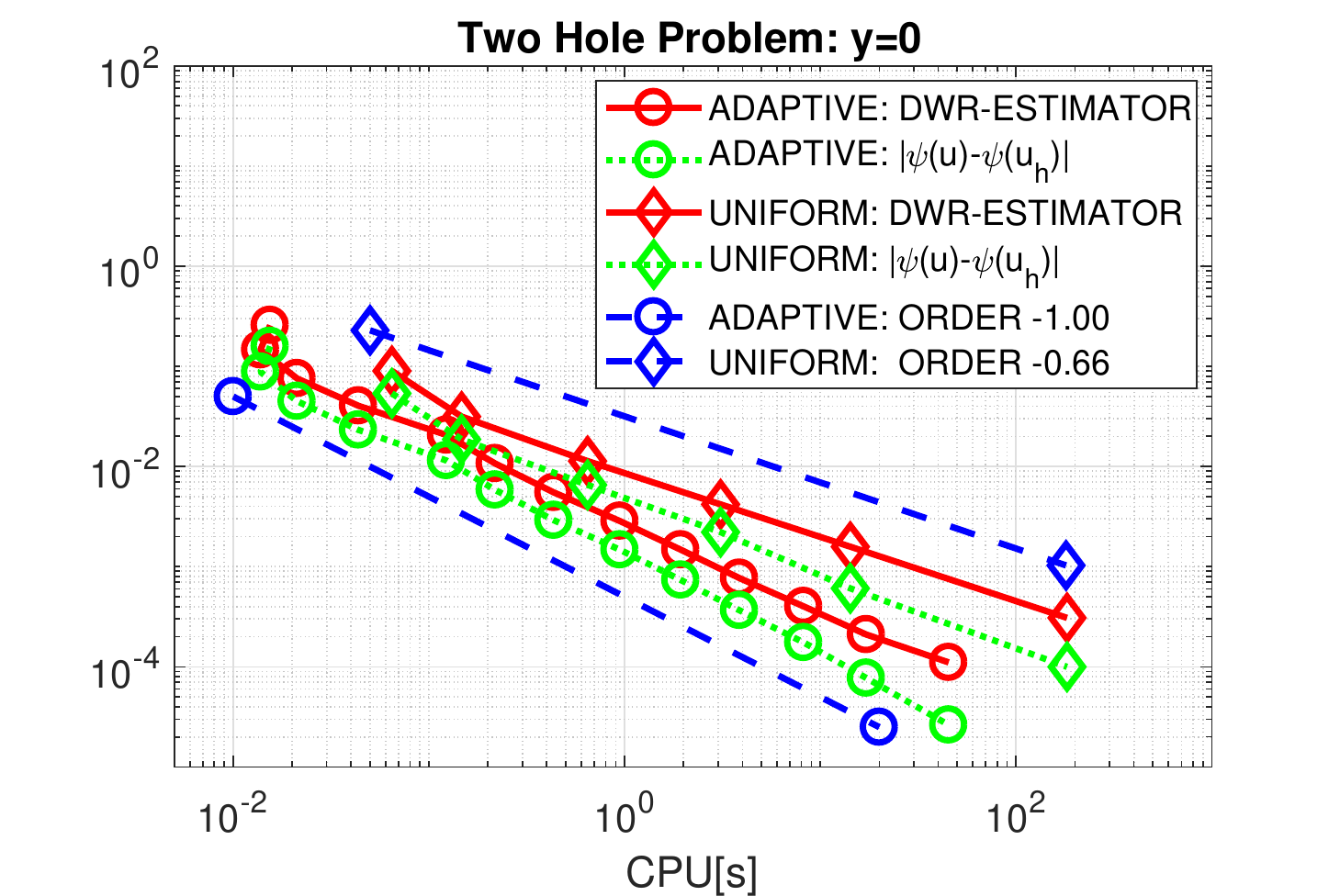}\includegraphics[width=0.55\textwidth]{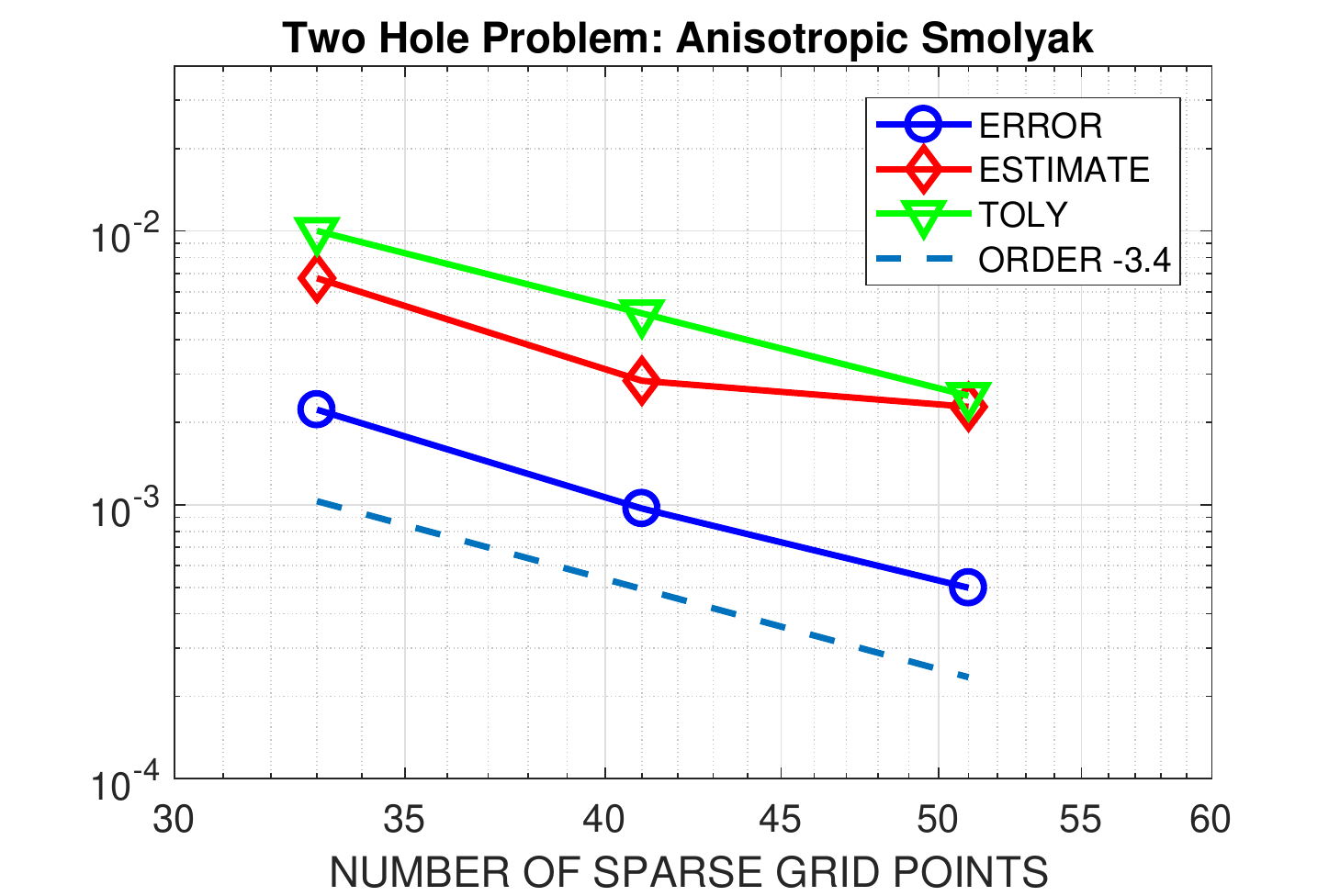}
\parbox{13cm}{
\caption{\small Two hole problem: History of error estimators and true errors
for $y=0\in\R^{16}$ using adaptive and uniform refinements
(left). A spatial tolerance of $\TolX=10^{-4}$ yields $3\,740\,240$ adaptive mesh points.
The error estimators and the true errors for anisotropic
Smolyak approximations of the expected value are shown only for low
tolerances  $\TolX=\TolY=10^{-2}/2^i,i=0,1,2$.
The numerically observed convergence orders for $\psi(u(y_0))$ in terms of CPU time
are close to $-1$ and $-0.66$ for adaptive and uniform refinement, respectively.
The numerically observed, averaged convergence order for $\E[\psi(u)]$
in terms of collocation points is $-3.4$ (right).}
\label{fig:TwoHoleEstAdaptUniform}
}
\end{figure}

All calculations start with an initial criss-cross structured mesh consisting of
$1920$ triangles, which are adjusted to the random holes. In Fig.~\ref{fig:TwoHoleSolMesh},
the numerical solution and its corresponding
adaptive mesh for a spatial tolerance $\TolX=10^{-1}$ and the random vector
$y\!=\!(0.5,0.5,-0.5,-0.5,-0.5,-0.5,1,-1,-1,-1,1,1,-1,1,-1,-1)^T$ are shown.
The adaptive algorithm refines the mesh at the eight reentrant corners due
to the fact that the exact solution contains a loss of regularity there. Exemplarily,
we study the convergence rates of adaptive and uniform refinements for the nominal
domain $D_0$, where $u(0)$ is contained in $H^{5/3}(D_0)$. While the correctly
adapted grids still recover the optimal order $-1$ for the approximation of
$\psi(u(0))$ in terms of CPU time, the use of a sequence of uniform meshes sees
the order drop down to the theoretical value $-0.66$, as can be seen in the left
part of Fig.~\ref{fig:TwoHoleEstAdaptUniform}. The DWR-estimators perform
quite well and deliver accurate upper bounds of the error. So we again set
$C_X=1$ in (\ref{qoi_cx}) and use $s^*\!=\!1$ for adaptive meshes and $s^*=1.6$
for uniform meshes in assumption (A1). Note that the latter choice
takes into account that due to the larger interior angle at the random holes for some
parameter values $y \not=0$, the regularity of the solution, and thus also the spatial convergence
rate, is slightly lower.

In order to estimate the parameter $\mu^*$, we apply the anisotropic
Smolyak algorithm and calculate samples at reasonable costs with tolerances
$\TolX=\TolY=10^{-2}/2^i,i=0,\ldots,3$. The value of $\E[\psi(u)]$
for $i=3$ is taken as the reference value, leading to an approximate convergence order of
$\mu^*=3.4$. From Fig. \ref{fig:TwoHoleEstAdaptUniform}, we detect that the
errors are significantly smaller than the prescribed tolerances, as it was also the
case in the first test example. Therefore, we again choose $C_Y=0.1$ and set
$\mu^*=3.4$ in assumption (A2). For later comparison, we compute
a reference value $\E[\psi(u)]=4.758057\ldots$
with $\TolX=10^{-4}$ and $\TolY=5\times10^{-4}$. The required
number of collocation points is $233$ and the
vector of maximum polynomial degree is
$(3\,2\,3\,3\,3\,3\,3\,3\,3\,3\,2\,2\,3\,2\,2\,2)^T$.

For the adaptive multilevel approach, we use again three levels 
and set the reduction factor to $q=0.2$. The spatial tolerances are then
computed from $\TolXk=\epsilon q^{k-2}/2$, $k=0,1,2$, with the six overall
accuracy requirements $\epsilon=10^{-2},5\times10^{-3},2.5\times10^{-3},10^{-3},5\times10^{-4},2.5\times10^{-4}$.
Following the algorithm given in Table~\ref{tab-algo-amlsc}, we first compute
$\TolXm1=\E[\psi(u_0)]$ with $\TolX=\TolY=\TolXzero$ and then derive the
stochastic tolerances $\TolYk$ from (\ref{opt-toly}) with
$C_Y=0.1$, $\mu^*=3.4$, and $s^*=1$.
The corresponding values are given in Table~\ref{tab-amlcol-two-hole}. For
the multilevel approach with uniform spatial meshes, we use $s^*=1.6$ to
determine the stochastic tolerances.

\begin{table}[t!]
\centering
\small
\begin{tabular}[t]{|l|c|c|c|c|c|c|}
\hline\rule{0mm}{2mm}\hspace{-0.1cm}
$\epsilon$ & $10^{-2}$ & $5\times10^{-3}$ & $2.5\times10^{-3}$ & $10^{-3}$ &
  $5\times10^{-4}$ & $2.5\times10^{-4}$\\
\hline\rule{0mm}{2mm}\hspace{-0.1cm}
$\TolXm1$ & $4.7125$ & $4.7346$ & $4.7465$ & $4.7525$ & $4.7554$ & $4.7566$ \\
\hline\rule{0mm}{2mm}\hspace{-0.1cm}
$\TolXzero$    & $1.3\times10^{-1}$ & $6.3\times10^{-2}$ & $3.1\times10^{-2}$ & $1.3\times10^{-2}$
  & $6.3\times10^{-3}$ & $3.1\times10^{-3}$ \\
$\TolXone$    & $2.5\times10^{-2}$ & $1.3\times10^{-2}$ & $6.3\times10^{-3}$ & $2.5\times10^{-3}$
  & $1.3\times10^{-3}$ & $6.3\times10^{-4}$ \\
$\TolXtwo$    & $5.0\times10^{-3}$ & $2.5\times10^{-3}$ & $1.3\times10^{-3}$ & $5.0\times10^{-4}$
  & $2.5\times10^{-4}$ & $1.3\times10^{-4}$ \\
\hline\rule{0mm}{2mm}\hspace{-0.1cm}
$\TolYtwo$    & $7.2\times10^{-3}$ & $4.1\times10^{-3}$ & $2.3\times10^{-3}$ & $1.1\times10^{-3}$
  & $6.2\times10^{-4}$ & $3.5\times10^{-4}$ \\
$\TolYone$    & $1.3\times10^{-2}$ & $6.1\times10^{-3}$ & $3.0\times10^{-3}$ & $1.1\times10^{-3}$
  & $5.5\times10^{-4}$ & $2.6\times10^{-4}$ \\
$\TolYzero$    & $3.0\times10^{-2}$ & $1.5\times10^{-2}$ & $7.2\times10^{-3}$ & $2.8\times10^{-3}$
  & $1.3\times10^{-3}$ & $6.4\times10^{-4}$ \\
\hline\rule{0mm}{2mm}\hspace{-0.1cm}
$M_0$ & $35$ & $51$ & $65$ & $131$ & $161$ & $321$ \\
$M_1$ & $33$ & $33$ & $33$ & $33$ & $33$ & $33$ \\
$M_2$ & $33$ & $33$ & $33$ & $33$ & $33$ & $33$ \\
\hline
\end{tabular}\\
\parbox{13cm}{
\caption{\small Adaptive multilevel collocation method: Sequences of
spatial and stochastic tolerances, $\TolXk$ and $\TolYk$, for
overall tolerances $\epsilon=10^{-2},\ldots,2.5\times10^{-4}$ and
certain approximations of $\TolXm1$ (above). Number of
collocation points taken by the anisotropic Smolyak algorithm
for $K\!=\!2$ (below).}
\label{tab-amlcol-two-hole}
}
\end{table}

\begin{figure}[t]
\centering
\includegraphics[width=0.8\textwidth]{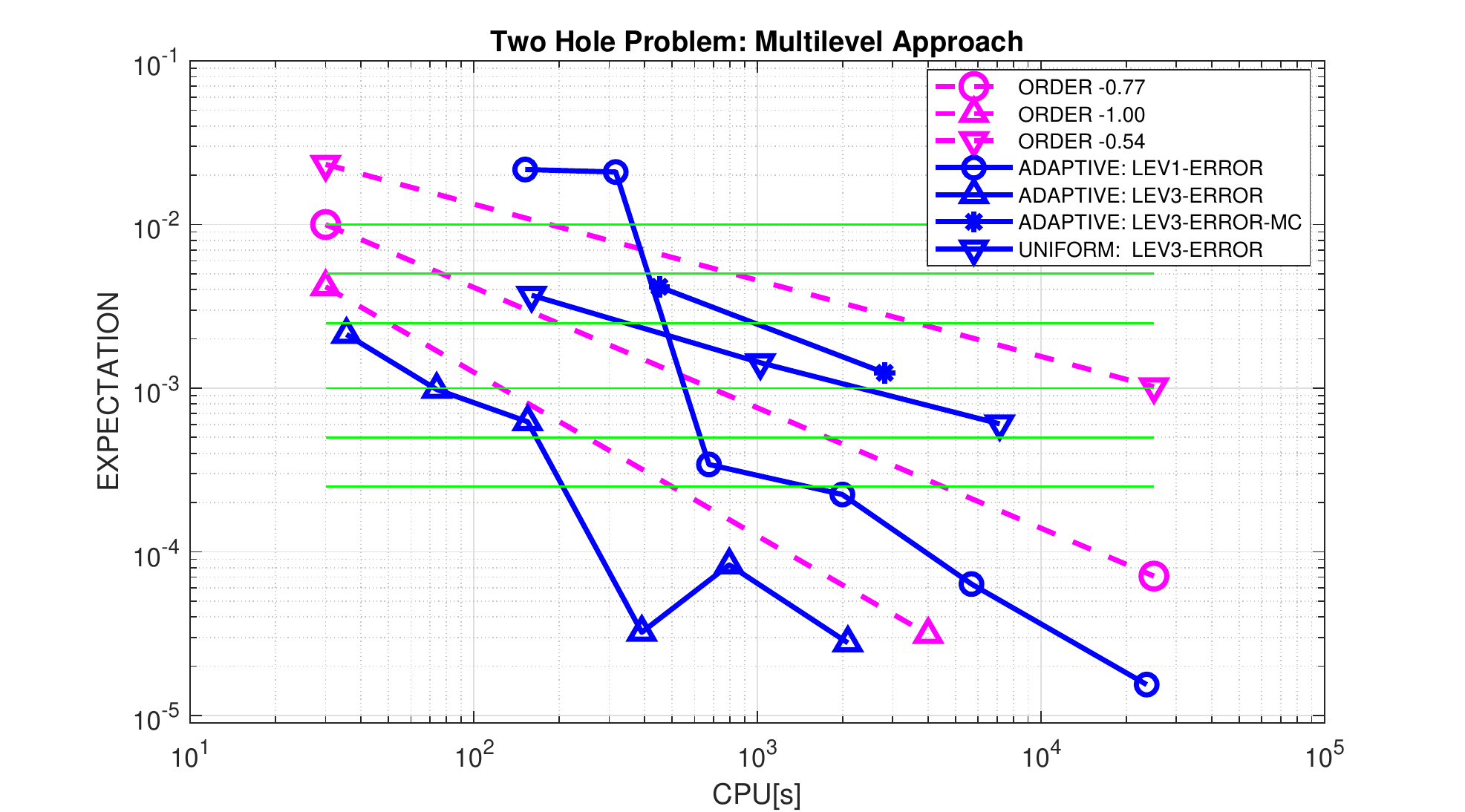}
\parbox{13cm}{
\caption{\small Two hole problem: Errors for the expected values $\E[\psi_{2}^{(\ML)}]$
and $\E[\psi_{2}^{(\SL)}]$ for the three-level (blue triangles)
and the one-level
(blue circles) approach with adaptive spatial meshes for
$\epsilon= 10^{-2},5\times10^{-3},2.5\times10^{-3},10^{-3},5\times10^{-4},2.5\times10^{-4}$ (green lines).
The orders of convergence predicted by Theorem \ref{th:qoi} are $-1$ and $-0.77$, respectively.
Results for $\epsilon\le 10^{-3}$ are not shown due to
excessive memory requirements ($>1.6\times10^{7}$ nodes). For comparison,
the results (adaptive: lev3-error-mc) of two runs of a three-level
adaptive multilevel Monte  Carlo method for $\epsilon=
10^{-2},5\times10^{-3}$ are also shown.}
\label{fig:TwoHoleMultiLevel}
}
\end{figure}

In Fig.~\ref{fig:TwoHoleMultiLevel}, errors for the expected value
$\E[\psi_2^{(\ML)}]$ versus computing time are plotted for the three-level
approach with adaptive and uniform spatial meshes. For comparison, we also
show results for the single-level approach with adaptive spatial meshes,
computed with $\TolX=\epsilon/2$ and $\TolY=\epsilon/(2C_Y)$. In
all cases, except the first two for the single-level method, the overall
tolerances are satisfied. The schemes work remarkably reliably.
However, the achieved accuracy is often significantly better than the prescribed
one -- never more than a factor $10$ though. One possible reason could be that cancellation effects are typically
overlooked if the overall error is split into two parts, which are then
individually controlled.

The multilevel method with adaptive spatial meshes outperforms the single-level
approach. The number of collocation points taken by the anisotropic
Smolyak algorithm are listed in Table~\ref{tab-amlcol-two-hole}. We observe
that the number of samples for the differences keeps constant, showing that
the increasing samples in the zeroth level always catch enough information
to eventually reach the tolerance. The corresponding numbers of
collocation points for the single-level method are $(1, 1, 33, 41, 51, 105)$.
Note that the computing time also includes the effort for the estimation
process, which is not visible in these numbers; see also the discussion for our
realization of the adaptive anisotropic Smolyak algorithm above. It is
also clear that the averaged constant $C_Y=0.1$ is too optimistic for the
first runs with $\epsilon=10^{-2},5\times10^{-3}$, which leads to the algorithm
taking only one collocation point.
The approximate orders of convergence for both methods,
$p_{\ML}=-1/s^*=-1$ and $p_{\SL}=-1/(s^*+1/\mu^*)\approx-0.77$, predict the
observed asymptotic rates for the computing times quite
well. However, the actual estimates can only serve as rough
indicators for the achieved accuracy.

The multilevel method with uniform
spatial meshes performs better than the single-level one for the first two
tolerances, but becomes quickly inefficient for higher tolerances. Due
to the larger value $s^*=1.6$, it needs significantly more samples for coarser meshes:
$M_0=(65,131,567)$ and $M_1=(33,33,105)$, compared to the multilevel
approach with adaptive spatial meshes, see Table~\ref{tab-amlcol-two-hole}.
$M_2$ remains the same. The observed convergence order $-0.54$ is close to
the predicted value $-1/s^*=-0.625$. Also for this example, it becomes
obvious that uniform meshes cannot compete with adaptive meshes for
higher tolerances.

This conclusion is also valid for tests with adaptive multilevel Monte Carlo.
Results for tolerances $\epsilon=10^{-2},5\times10^{-3}$ are shown in
Fig.~\ref{fig:TwoHoleMultiLevel} for three levels. The numbers of optimized
samples are $M_0=(4140,17013)$, $M_1=(48,108)$, $M_2=(5,6)$, respectively.
All values are calculated from an average over 5 independent realizations.
Although the variance reduction is quite high, leading to surprisingly small
numbers $M_2$, the fast increasing numbers for $M_0$ are still too challenging,
especially for higher tolerances.

\section{Concluding remarks}
Adaptive methods have the potential to drastically reduce the number of
degrees of freedom and to realise optimal complexity in terms of computing time,
even in cases where the exact solutions have spatial singularities.
As well as providing a posteriori error estimates, adaptive methods usually outperform
methods based on uniform mesh refinement. When coupled with multilevel
stochastic algorithms, they can increase the computational efficiency in the
sampling process when the stochastic dimension increases and hence provide
a general tool to further delay the curse of dimensionality.

State-of-the-art
adaptive methods implemented in the open-source \matlab\ packages {\sl p1afem}
and {\sl Sparse Grid Kit} have been used in off-the-shelf fashion in this study.
The numerical results for the adaptive multilevel collocation method demonstrate
the reliability of the error control and a  significant decrease
in complexity compared to  uniform spatial refinement methods and single-level stochastic
sampling methods. While these
advantages are already known for the  methods involved, this study shows that
the adaptive components can be combined to give an efficient algorithm for solving PDEs
with random data to  arbitrary levels of accuracy.

\section*{Appendix.}
The pathwise solution of  the one peak problem is given by
$$
u({\vfld x},\bfy) = \exp ( - \beta  \{  \alpha(y_1) (x_1 -y_1)^2 + (x_2 -y_2)^2 \} ),  \eqno{(A.1)}
$$
with $\alpha(y_1)= 18 y_1 + 11/2$.
The quantity of interest  is given by
$$
\Bbb{E}\left[\phi(u)\right]= \int_\Gamma \int_D u^2(\vfld x,\bfy) \rho(y) \, d\vfld x \, d\bfy .  \eqno{(A.2)}
$$
The specific choice $\beta=50$   taken in the numerical experiments leads to two simplifications.
\begin{itemize}
\item[(S1)]
The Dirichlet boundary condition ($u$ satisfying (A.1) on $\partial D$) may be replaced
 without significant loss of accuracy by the numerical approximation
$$ u_h=0  \hbox{ on } \partial D.  \eqno{(A.3)}$$
{\bf Justification}: By direct computation, the maximum value of $u$ on the boundary
occurs for $\alpha=1$, when the standard deviation $\sigma$ of the Gaussian peak is
minimised (thus $\sigma_* = 1/10$ and $\bfy=(-1/4,y_2)$ for all values $-1/4\leq y_2\leq 1/4$).
For all  such values of $y_2$ the  nearest point  on the boundary is ${\vfld x} = (-1,y_2)$, and the
shortest distance to the boundary is $d=3/4$ so that the maximum value of
$u$ on the boundary is given by $ u_{\max} = \exp\big(- \frac{9\beta}{16}\big) \approx 6\cdot 10^{-13}$.
\item[(S2)]
Thus, we may readily compute a reference value  (accurate to more than 10 digits):
$$ \Bbb{E}\left[\phi(u)\right]\approx  Q := {1\over 9} \cdot (\sqrt{10 }-1)\cdot {\pi\over \beta}
= 0.015 095 545 \ldots$$

 {\bf Justification}:
 the integral in (A.2) is separable, thus
\begin{align*}
 I &=  \int_\Gamma \int_D u^2(\vfld x,\bfy) \rho(y)\, d\vfld x \, d\bfy \\
 &=
   \int_\Gamma \int_D e^{-2\beta \{ \alpha (x_1-y_1)^2  + (x_2-y_2)^2  \}  } \rho(y)\,  d\vfld x \, d\bfy \\
 &=
   \int_{-{1\over 4}}^{1\over 4} \hat{\rho}(y_1) \int_{-1}^{1} e^{-2\beta\alpha (x_1-y_1)^2 } \, dx_1\, dy_1 \cdot
   \int_{-1/4}^{1/4} \hat{\rho}(y_2) \int_{-1}^{1} e^{-2\beta  (x_2-y_2)^2 } \, dx_2\, dy_2 .
  \end{align*}
The key point is that the inner integrals in the above expression are  Gaussians with small variance.
In both cases the integrand is close to zero at the end points $-1$ and $1$. Thus, in both
cases the  integral can be extended by zero to the range $(-\infty,\infty)$ (this could be made
rigorous using probabilistic arguments) to give the following  reference value approximation
\begin{align*}
 Q &=
   \int_{-{1\over 4}}^{1\over 4} 2 \int_{-\infty}^{\infty} e^{-2\beta\alpha (x_1-y_1)^2 } \, dx_1\, dy_1 \cdot
   \int_{-1/4}^{1/4} 2 \int_{-\infty}^{\infty} e^{-2\beta  (x_2-y_2)^2 } \, dx_2\, dy_2  \\
   &=   4 \int_{-{1\over 4}}^{1\over 4} \sqrt{\pi \over 2 \alpha \beta}\, dy_1 \cdot
   \int_{-1/4}^{1/4} \sqrt{\pi \over  2 \beta} \, dy_2  \\
   &=  {2 \pi\over \beta} \int_{-1/4}^{1/4}  \, dy_2   \cdot
    \int_{-{1\over 4}}^{1\over 4}  \alpha(y_1)^{-1/2} \, dy_1
   =  { \pi\over \beta}  \cdot {1\over 18}\,
   \underbrace{\int_{1}^{10}  \alpha^{-1/2}  \,  d\alpha}_{2(\sqrt{10}-1)}  .
 \end{align*}
 \end{itemize}

\medskip\noindent
{\bf Acknowledgements.}
The first author is supported by the Deutsche Forschungsgemeinschaft
(DFG, German Research Foundation) within the collaborative research center
TRR154 {\em ``Mathematical modeling, simulation and optimisation using
the example of gas networks''} (Project-ID 239904186, TRR154/2-2018, TP B01),
the Graduate School of Excellence Computational Engineering (DFG GSC233),
and the Graduate School of Excellence Energy Science and
Engineering (DFG GSC1070). The second author was supported by the
UK Engineering and Physical Sciences Research Council (EPSRC) under
grant number EP/K031368/1. The third author was funded by a
Romberg visiting scholarship at the University of Heidelberg in 2019.
The authors would also like to thank the Isaac Newton Institute for
Mathematical Sciences, Cambridge, for support and hospitality during
the programme {\em ``Uncertainty quantification for complex systems: theory
and applications''} (Jan--Jul 2018), where  the work on this
paper was initiated.

\bibliographystyle{plain}
\bibliography{bibmluq}

\end{document}